\documentclass[MSNbibl,number,citesort,seceqn,dvips]{arxbj}
\usepackage{upgreek,mathbh}
\usepackage{graphicx}

\aid{0}
\volume{19}
\issue{5A}
\pubyear{2013}
\firstpage{1576}
\lastpage{1611}
\doi{10.3150/12-BEJ422} 

\makeatletter

\define@key{bid}{pmid}{}

\newproclaim{assumption}{Assumption}
\newtheorem{theo}{Theorem}[section]
\newremark{rem}[theo]{Remark}
\newtheorem{coro}[theo]{Corollary}
\newtheorem{prop}[theo]{Proposition}
\newproclaim{defin}[theo]{Definition}
\newtheorem{lem}{Lemma}

\renewcommand{\phi}{\varphi}
\renewcommand{\epsilon}{\varepsilon}
\newcommand{\E}{\mathbb{E}}
\renewcommand{\P}{\mathbf{P}}
\newcommand{\Ex}{\mathbb{E}}
\newcommand{\var}{\operatorname{\mathbb{V}ar}}
\newcommand{\argmin}{\operatorname{argmin}}
\newcommand{\argmax}{\operatorname{argmax}}
\newcommand{\cc}[1]{\overline{#1}} 
\newcommand{\1}{\mathbh{1}}
\newcommand{\diag}{\operatorname{diag}}
\newcommand{\vek}[1]{\underline{#1}}
\newcommand{\C}{{\mathbb C}}
\newcommand{\N}{{\mathbb N}}
\newcommand{\R}{{\mathbb R}}
\newcommand{\Z}{{\mathbb Z}}
\newcommand{\cB}{{\mathcal B}}
\newcommand{\cD}{{\mathcal D}}
\newcommand{\cF}{{\mathcal F}}
\newcommand{\cR}{{\mathcal R}}
\newcommand{\cS}{{\mathcal S}}
\newcommand{\td}{\widetilde{d}}
\newcommand{\tf}{\widetilde{f}}
\newcommand{\tk}{\widetilde{k}}
\newcommand{\tw}{\widetilde{w}}
\newcommand{\hf}{\widehat{f}}
\newcommand{\hg}{\widehat{g}}
\newcommand{\supt}[1]{\sup_{t\in\cB_{#1}}}
\newcommand{\Mnm}{M_{n,m}}
\newcommand{\pen}{\operatorname{pen}}
\newcommand{\hpen}{\widehat{\pen}}
\newcommand{\hPhi}{\widehat\Phi}
\newcommand{\tPhi}{\widetilde\Phi}
\newcommand{\ct}{\Upsilon}
\newcommand{\whk}{\widehat k}
\newcommand{\kstar}{k_n^*}
\newcommand{\kstarm}{k_m^*}
\newcommand{\Cy}{N}
\newcommand{\Ce}{M}
\newcommand{\hCy}{\widehat N}
\newcommand{\hCe}{\widehat M}
\newcommand{\xdfw}{\gamma}
\newcommand{\xdfr}{r}
\newcommand{\Edfw}{\Lambda}
\newcommand{\edfr}{d}
\newcommand{\xdf}{f}
\newcommand{\ydf}{g}
\newcommand{\edf}{\varphi}
\newcommand{\hxdf}{\widehat\xdf}
\newcommand{\hydf}{\widehat\ydf}
\newcommand{\hedf}{\widehat\edf}
\newcommand{\txdf}{\widetilde\xdf}
\newcommand{\ef}{e}
\newcommand{\bw}{\xdfw}
\newcommand{\br}{\xdfr}
\renewcommand{\td}{d}
\renewcommand{\tw}{\lambda}
\newcommand{\hw}{\omega}
\newcommand{\hDelta}{\widehat{\Delta}}
\newcommand{\hk}{\whk}
\renewcommand{\Mnm}{(\Cy_n\wedge\Ce_m)}
\newcommand{\Mmn}{\Ce_{m_n}}
\newcommand{\Fgr}{\cF_\xdfw^\xdfr}
\newcommand{\Eld}{{\mathcal E}_\lambda^\edfr}
\newcommand{\hfedf}{\widehat{[\edf]}}
\newcommand{\hfydf}{\widehat{[\ydf]}}
\newcommand{\tCe}{\widetilde{\Ce}}
\newcommand{\mod}{\operatorname{mod}}
\newcommand{\eqref}[1]{(\ref{#1})}
\makeatother

\begin{document}
\begin{frontmatter}

\title{Adaptive circular deconvolution by model selection
under unknown error distribution}
\runtitle{Adaptive circular deconvolution}

\begin{aug}
\author{\fnms{Jan} \snm{Johannes}\corref{}\thanksref{e1}\ead[label=e1,mark]{jan.johannes@uclouvain.be}} \and
\author{\fnms{Maik} \snm{Schwarz}\thanksref{e2}\ead[label=e2,mark]{maik.schwarz@uclouvain.be}}
\runauthor{J. Johannes and M. Schwarz} 
\address{Institut de statistique, biostatistique et sciences actuarielles,
Voie du Roman Pays 20, Bo\^ite L1.04.01, 1348 Louvain-la-Neuve,
Belgium. \printead{e1},\\ \printead*{e2}}
\end{aug}

\received{\smonth{1} \syear{2011}}
\revised{\smonth{8} \syear{2011}}

%
\begin{abstract}
We consider a circular deconvolution problem, in which the density
$\xdf$ of a circular
random variable $X$ must be estimated nonparametrically based on an
i.i.d. sample from a noisy observation $Y$ of $X$. The additive
measurement error is supposed to be independent
of $X$. The objective of this work was to construct a fully data-driven
estimation procedure when the error density
$\phi$ is unknown. We assume that in addition to the i.i.d.
sample from $Y$, we have at our disposal an additional
i.i.d. sample drawn independently from the error
distribution. We first develop a minimax theory in terms of both sample
sizes. We
propose an orthogonal series estimator attaining the minimax rates but
requiring optimal choice of a dimension parameter depending on
certain characteristics of $f$ and $\phi$, which are not known in
practice. The main issue addressed in this work is the adaptive choice
of this dimension parameter using a model selection approach. In a
first step, we develop a penalized minimum contrast estimator
assuming that the error density is known.
We show that this partially adaptive estimator can attain the lower risk
bound up to a constant in both sample sizes $n$ and $m$.
Finally, by randomizing the penalty and the
collection of models, we modify the estimator such that it no longer
requires any previous knowledge of the error distribution. Even
when dispensing with any hypotheses on $\phi$, this fully data-driven
estimator still preserves minimax optimality in almost the same cases
as the partially adaptive estimator.
We illustrate our results by computing minimal rates under
classical smoothness assumptions.
\end{abstract}

%

\end{frontmatter}

\section{Introduction}\label{sec:intro}
\begin{figure}

\includegraphics{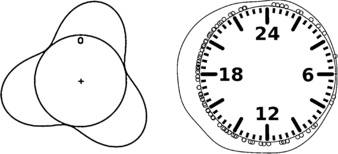}

\caption{A trimodal circular density and a density estimator from
periodic data.}\label{fig:density}
\end{figure}

This work deals with the estimation of circular probability densities
from noisy observations. ``Circular'' means that the
observations are points on the circle. Such models arise in
numerous and various fields of application. Data with temporal
structure are most naturally represented in this way; for example,
times of day when events of interest occur such as requests in a
computer network, financial transactions, or gun
crimes, can be represented as points on a clock face (Gill and
Hangartner \cite{GH:10}), as illustrated in Figure~\ref{fig:density}.
Replacing the clock face by a compass rose,
directional data also can be treated in the circular
setting. Curray \cite{Cur:56} considered the analysis of directional
data in
the context of geological research.
Cochran, Mouritsen and Wikelski \cite{CMW:04} investigated
migrating birds' navigation abilities using circular
data.

%
%

The applications of circular data are not restricted to a
spatiotemporal context. Gill and Hangartner \cite{GH:10} provided an
overview of circular
data in political science, where they can be used to, for example, model
political preferences, which are not of a temporal or a spatial
nature. For a more detailed discussion of the specifics of
circular data, see Mardia \cite{Mardia1972}. Numerous circular data
sets and examples of their statistical analysis have been provided by
Fisher \cite{Fisher1993}.

Let $X$ be the circular random variable whose density $\xdf$ we are
interested in and let $\epsilon$ be an independent additive circular error
with unknown density $\edf$. Denote by $Y$ the contaminated observation
and by $\ydf$ its density. Throughout this work, we
identify the circle with the unit interval $[0,1)$ for
notational convenience. Thus, $X$ and $\epsilon$
take their values in $[0,1)$. Let $\lfloor\cdot\rfloor$ be the floor
function. Taking into account the circular nature of the data, the
model can be written as $Y=X+\epsilon-\lfloor X+\epsilon\rfloor$ or,
equivalently,
$Y=X+\epsilon\mod\ [0,1)$.
We then have
\[
\ydf(y) = (\xdf\ast\edf) (y):= \int_{[0,1)} \xdf\bigl((y-s) -
\lfloor y-s\rfloor\bigr) \edf(s) \,\mathrm{d}s,\qquad y\in[0,1),
\]
such that $\ast$ denotes circular convolution. Therefore, the
estimation of $\xdf$ is called a circular deconvolution problem. Let
$L^2:=L^2([0,1))$ be the Hilbert space of square-integrable
complex-valued functions defined on $[0,1)$ endowed with the usual
inner product $\langle f,g\rangle=\int_{[0,1)} f(x) \overline{
g(x)}\,\mathrm{d}x$,
where $\overline{ g(x)}$ denotes the complex conjugate of ${
g(x)}$. In this work, we suppose that $\xdf$ and $\edf$, and hence also
$\ydf$, belong to the subset $\cD$ of all densities in $L^2$.
Consequently, they admit representations as discrete Fourier series
with respect to the exponential basis, $\{e_j\}_{j\in\Z}$, of $L^2$,
where $ e_j(x): = \exp(-\mathrm{i}2\uppi j x)$ for $x\in[0,1)$ and
$j\in\Z$. Given $p\in\cD$ and $j\in\Z$, let $[p]_j:=\langle
p,e_j\rangle$
be the $j$th Fourier coefficient of $p$. In particular, $[p]_0 =1$.
The key to the analysis of the circular deconvolution problem is the
convolution theorem, which states that $g=f\ast\phi$ if and only if
$[\ydf]_j = [\xdf]_j [\edf]_j$ for all
$j\in\Z$. Therefore, as long as $[\edf]_j\ne0$ for all $j\in
\Z$,
which we assume from here on, we have
%
%
\begin{equation}
\label{eq:14} \xdf= 1 + \sum_{\vert j\vert>0} \frac{[\ydf
]_j}{[\edf]_j}
e_j \qquad\mbox{with } [\ydf]_j=\Ex
e_j(-Y) \mbox{ and } [\edf]_j=\Ex e_j(-
\epsilon)\ \forall j\in\Z.
\end{equation}

Note that an analogous representation holds in the case of
deconvolution on the real line when the $X$-density is compactly
supported but the error term $\epsilon$, and hence $Y$, take their
values in $\R$. In this situation, the deconvolution density still
admits a discrete representation as in \eqref{eq:14}, but involving
the characteristic functions of $\edf$ and $\ydf$ rather than their
discrete Fourier coefficients. There is a vast literature on
deconvolution on the real line, with or without compactly supported
deconvolution density. In the case where the error density is fully
known, a very popular approach based on kernel methods has been
considered by,
among many others, Carroll and Hall \cite{CarrollHall88}, Devroye
\cite{Devroye89},
Fan \cite{Fan91,Fan92},
Stefanski \cite{Stefanski90}, Zhang \cite{Zhang90}, Goldenshluger
\cite{Goldenshluger99,Goldenshluger00}, and Kim and Koo \cite{KimKoo02}.
Mendelsohn and Rice \cite{MendelsohnRice82} and Koo and Park \cite{KooPark96},
for example,
studied spline-based methods, whereas
Pensky and Vidakovic \cite{PenskyVidakovic99}, Fan and Koo \cite
{FanKoo02}, and
Bigot and Van~Bellegem \cite{BigotVanBellegem2009}, used wavelet decomposition.
Situations
with only partial knowledge of the error density have been considered
as well
(e.g., Butucea and Matias \cite{ButuceaMatias2005}, Meister \cite
{Meister2006},
Schwarz and Van~Bellegem \cite{SchwarzVanBellegem2009}). Consistent
deconvolution without previous
knowledge of the error distribution is also possible in the case of
panel data (e.g., Horowitz and Markatou \cite{HorowitzMarkatou1996},
Hall and Yao \cite{HallYao2003}, or
Neumann \cite{Neumann2007}) or by assuming an additional sample from
the error
distribution (e.g., Diggle and Hall \cite{DiggleHall1993}, Neumann
\cite{Neumann1997},
Johannes \cite{Jo07}, or Comte and Lacour \cite{ComteLacour2009}).
For a broader overview on
deconvolution problems, see the monograph of
Meister \cite{Meister2009}.

We now return to the circular case. In this paper, we assume that we
do not know the density $g=\xdf\ast\edf$ of the contaminated
observations or the error density $\edf$, but we have at our
disposal two independent samples of i.i.d. random variables
%
%
\begin{equation}
\label{eq:obs} Y_k \sim g \qquad(k=1,\ldots,n) \quad\mbox{and}
\quad\epsilon_k \sim\edf\qquad(k=1,\ldots,m)
\end{equation}
of size $n\in\N$ and
$m\in\N$, respectively. Our aim is to establish a fully data-driven
estimation procedure for the deconvolution density $\xdf$ that attains
optimal convergence rates in a minimax sense. More precisely,
given classes $\Fgr$ and $\Eld$ (defined below) of deconvolution and
error densities, respectively, we measure the accuracy of an
estimator $\txdf$ of $\xdf$ by the maximal weighted risk
$\sup_{\xdf\in\Fgr}\sup_{\edf\in\Eld}\Ex\Vert\txdf-\xdf
\Vert_\hw^2$
defined with respect to some weighted norm $\Vert\cdot\Vert^2_\hw:=
\sum_{j\in\Z}\hw_j\vert[\cdot]_j\vert^2$, where $\hw:= (\hw
_j)_{j\in\Z}$
is a strictly positive sequence of weights. This allows us to quantify
the estimation accuracy in terms of the mean integrated squared error
(MISE) not only of $\xdf$ itself, but also of its derivatives, for
example. It is well known that even in case of a known error density,
the maximal risk in terms of the MISE in the circular deconvolution
problem is essentially determined by the asymptotic behavior of the
sequences of Fourier coefficients $([\xdf])_{j\in\Z}$ and
$([\edf])_{j\in\Z}$ of the deconvolution density and the error
density, respectively. For a fixed deconvolution density $\xdf$, a
faster decay of the $\epsilon$-density's Fourier coefficients
$([\edf])_{j\in\Z}$ results in a slower optimal rate of
convergence. For example, in the standard context of an ordinary smooth
deconvolution density, when $([\xdf])_{j\in\Z}$
decays polynomially, logarithmic rates of convergence appear when the
error density is super smooth, that is, $([\edf])_{j\in\Z}$ has
exponential decay. Efromovich \cite{Efromovich97} treated this
special case exclusively. However, this situation and many others
are covered by the density classes
\begin{eqnarray*}
\cF_\xdfw^\xdfr&:=& \biggl\{ p\in\cD\dvt\sum
_{j\in\Z} \xdfw_j \bigl\vert[p]_j
\bigr\vert^2 =:\Vert p\Vert_\xdfw^2\leq\br\biggr\}
\quad\mbox{and}
\\
\Eld&:=& \biggl\{ p\in\cD\dvt1/\edfr\leq\frac
{\vert[p]_j\vert^2}{\lambda_j}\leq\edfr\ \forall j
\in\Z\biggr\},
\end{eqnarray*}
where $\xdfr,\edfr\geq1$ and the positive weight sequences $\xdfw:=
(\xdfw_j)_{j\in\Z}$ and $\lambda:= (\lambda_j)_{j\in\Z}$ specify the
asymptotic behavior of the respective sequence of Fourier
coefficients. In Section~\ref{sec:minimax}, we present a lower bound of
the maximal weighted risk that is determined essentially by the
sequences $\xdfw$, $\lambda$, and $\hw$. This lower bound is
composed of
two main terms, each of which depends on the size of one sample but
not of the other sample. Let us define an orthogonal series estimator
by replacing the unknown Fourier coefficients in (\ref{eq:14}) by
empirical counterparts, that is,
%
%
\begin{eqnarray}
\label{def:est} &&\hxdf_k:=1+\sum_{0<\vert j\vert\leq k}
\frac{\widehat{[\ydf]}_j}{\widehat{[\edf]}_j}\mathbh{1} { \bigl\{
\bigl\vert
\widehat{[\edf]}_j
\bigr\vert^2 \geq1/m \bigr\}} \ef_j
\nonumber
\\[-8pt]
\\[-8pt]
&&\quad\mbox{with } \widehat{[\ydf]}_j := \frac{1}{n}\sum
_{i=1}^n e_j(-Y_i) \mbox{ and } \widehat{[\edf]}_j := \frac
{1}{m}\sum
_{i=1}^m e_j(-\epsilon_i).
\nonumber
\end{eqnarray}
For each $j$, we introduce a threshold for the estimated coefficient
$\widehat{[\edf]}_j$ that corresponds, in accordance with Neumann
\cite
{Neumann1997}, to
the rate at which $[\edf]_j$ can be estimated. Again, things work
out analogously to deconvolution on the real line, where we need only
replace the empirical Fourier coefficients with the corresponding values
of the empirical characteristic functions. Similar estimators have
been studied by, for example, Neumann \cite{Neumann1997} on the real
line and by
Efromovich \cite{Efromovich97} in the circular case.

We show below that the estimator $\hxdf_k$ attains the lower
bound and thus is minimax optimal. By comparing the minimax rates in
the cases of known and unknown error density, we can
characterize the influence of the estimation of the error
density on the quality of the estimation. In particular, depending on
the $Y$ sample size $n$, we can determine the minimal
$\epsilon$ sample size $m_n$ needed to attain the same upper risk
bound as in the case of a known error density, up to a
constant. Interestingly, the required sample size, $m_n$, is far smaller
than $n$ in a wide range of situations. For example, in the super
smooth case, it is sufficient that the size of the $\epsilon$ sample
be a polynomial in $n$, that is, $m_n=n^r$ for any $r>0$.

Of course, minimax optimality can be achieved only if the
dimension parameter $k$ is chosen in an optimal way. In general, this
optimal choice of $k$ depends on, among other things, the sequences
$\xdfw$
and $\lambda$. However, in the special case where the error density is
known to be super smooth and the deconvolution density is ordinary
smooth, the optimal dimension parameter depends only on $\lambda$ and
not on $\xdfw$. Thus, the estimator is automatically adaptive with
respect to $\xdfw$ under the optimal choice of $k$. In this situation,
Efromovich \cite{Efromovich97} provided an estimator that is also
adaptive with
respect to the super smooth error density. In contrast,
Cavalier and Hengartner~\cite{CavalierHengartner2005}, deriving oracle
inequalities in an
indirect regression problem based on a circular convolution
contaminated by Gaussian white noise, treated only the ordinary smooth
case. As in our setting, their observation scheme involves two
independent samples. Of note, application of these
estimators requires knowledge of whether the error density
is ordinary or super smooth. In this work, we provide a unified
estimation procedure that can attain minimax rates in both cases, being adaptive over a class including both
ordinary and super smooth error densities. This fully adaptive method
of choosing the parameter $k$ depends only on the observations, not
on characteristics of either $\xdf$ or $\edf$. Our main result
is that for this automatic choice $\whk$, the
estimator $\hxdf_{\whk}$ attains the lower bound up to a constant, and
thus is minimax-optimal, over a wide range of sequences $\xdfw$
and $\lambda$, covering in particular both ordinary and super smooth
error densities. A similar result was recently derived in the
context of a functional linear regression model by
Comte and Johannes~\cite{CJ:10}.

Regarding the two sample sizes, the assumption of Cavalier and Hengartner
\cite{CavalierHengartner2005} on the respective noise levels can be
translated to our model by stating that the $\epsilon$ sample size $m$
is at least as large as the $Y$ sample size $n$. This assumption was
also made by Efromovich \cite{Efromovich97}. Also note that in the functional
linear regression model, only one sample size, $n$, occurs (Comte and
Johannes \cite
{CJ:10}); however, as mentioned earlier,
without changing the minimax rates, the $\epsilon$-sample size can be
reduced to $m_n$, which can be much smaller than $n$. This is a
desirable property, given that the observation of the additional sample from
$\epsilon$ may be expensive in practice. Nevertheless, the minimal
choice of $m$ depends on, among other things, the sequences $\xdfw$
and $\lambda$ and thus is unknown in general. Despite the eventual
deterioration of the minimax
rate resulting from choosing the sample size $m$
smaller than $n$, the proposed estimator still attains this rate in
many cases; that is, no price, in terms of convergence rate, is paid
for adaptivity.

The adaptive choice of $k$ is motivated by the general model selection
strategy developed by Barron, Birg{\'e} and Massart \cite
{BarronBirgeMassart1999}. Concretely,
following Comte and Taupin \cite{ComteTaupin}, who treated the case of
a known error
density only, $\hk$ is the minimizer\footnote{For a sequence $a_n$
attaining a minimum on $N\subseteq\N$, let
$\argmin_{n\in N} a_n := \min\{n\in N \vert a_n \leq
a_k\ \forall k\in N\}$.} of a penalized contrast
\[
\whk:= \mathop{\argmin}_{1\leq k\leq K} \bigl[-\Vert\hxdf_{k}
\Vert^2_\hw+ \pen(k) \bigr].
\]
Note that we can compute $\Vert\hxdf_{k}\Vert^2_\hw=1+\sum
_{0<\vert j\vert\leq
k}\hw_j{\vert[\hydf]_j\vert^2}{\vert[\hedf]_j\vert^{-2}}\1\{
\vert[\hedf]_j\vert^2\geq
1/m\}$. As in case of a known error density, it turns out
that both the penalty function $\pen(\cdot)$ and the upper bound $K$
needed for the correct choice of $k$ depend on a
characteristic of the error density, which is now unknown. This
quantity is often referred to as the
degree of ill-posedness of the underlying inverse problem. Therefore,
as an intermediate step, we allow the penalty function $\pen(\cdot)$
and the upper bound
$K$ to depend on the error density. We then show
an upper risk bound for the resulting partially adaptive
estimator. We prove that over a wide range of sequences
$\xdfw$, this choice of $k$ yields the same upper
risk bound as the optimal choice, up to a constant. Finally, we choose
$k$ fully adaptively by replacing $\pen(\cdot)$ and
$K$ by their empirical versions, which depend only on the data. As in
the case of known
degree of ill-posedness, we show an upper risk bound for the now fully
adaptive estimator.

Let us return briefly to deconvolution on the real line with compactly
supported $X$ density. We note that in this situation, the adaptive
choice of $k$ can be performed in the same way. Moreover, the upper
risk bounds remain valid, and the adaptive estimator is minimax
optimal over a wide range of cases. In fact, the circular structure of
the model is exploited only in the proof of the lower bound and to
guarantee the existence of the discrete representation in
\eqref{eq:14}, which still holds in case of a compactly supported
deconvolution density.

This paper is organized as follows. In the next section, we develop
the minimax theory for the circular deconvolution model with respect to the
weighted norms introduced above and compute the rates which we can
obtain in different configurations for the weight
sequences. We devote the final section to constructing the adaptive
estimator and show an upper risk bound. We illustrate our results with
example configurations
considered in Section~\ref{sec:minimax}. All proofs are deferred to
the \hyperref[app]{Appendix}.

\section{Minimax optimal estimation}\label{sec:minimax}
In this section, we develop the minimax theory for estimating a
circular deconvolution density under unknown error density when two
independent samples from $Y$ and $\epsilon$, of size $n$ and $m$,
respectively, are available. We derive a lower bound depending on both sample
sizes and show that the orthogonal series estimator
$\hxdf_k$ defined in \eqref{def:est} attains this lower bound up to a
constant if $k$ is chosen in an appropriate way. All results in this
paper are derived under the following minimal regularity conditions:
\setcounter{assumption}{0}
\renewcommand{\theassumption}{A\arabic{assumption}}
\begin{assumption}\label{ass:minreg} Let $(\xdfw_j)_{j\in\Z}$,
$(\hw_j)_{j\in\Z}$ and $(\lambda_j)_{j\in\Z}$ be strictly
positive {symmetric} sequences of weights with
$\xdfw_0=\hw_0=\hw_1=\lambda_0 = \lambda_1=1$ such that $(\hw
_n/\bw_n)_{n\in\N}$ and $(\lambda_n)_{n\in\N}$ are
nonincreasing, respectively with $\Lambda:= \sum_{j\in\Z}\lambda_j
< \infty$.
\end{assumption}

Here and subsequently, we refer to any sequence $(a_n)_{n\in\Z}$
as a whole by omitting its index as in, for example ``the sequence
$a$''. We define arithmetic operations on sequences
element-wise.
Furthermore, we denote by $C$ universal numerical constants and
by $C(\cdot)$ constants depending only on the arguments. In both
cases, the values of the constants may change from line to line.
Moreover, we write
$a_n\lesssim b_n$ when $a_n\leq C
b_n$ for all sufficiently large $n\in\N$, and $a_n\sim b_n$ when
$a_n\lesssim b_n$ and $b_n\lesssim a_n$ simultaneously.

\subsection*{Lower bounds} The next assertion provides a lower bound in
the case of a known error density, which obviously will depend on the
size of the $Y$ sample only. Of course, this lower bound is still valid in
the case of an unknown error density.

\begin{theo}\label{theo:lower:n} Assume an i.i.d. $Y$ sample of
size $n$. Consider sequences $\hw$, $\xdfw$, and
$\lambda$ satisfying Assumption~\textup{\ref{ass:minreg}} such that
$\sum_{j\in\Z} \xdfw_j^{-1}=\Gamma<\infty$ and $\edf\in
\Eld$ for some $d\geq1$. Define, for all $n\geq1$,
%
%
\begin{eqnarray}
\label{def:psin} %
\kstar&:=&\kstar(\xdfw,\lambda,\hw):=
\mathop{\argmin}_{k\in\N} \biggl[\max\biggl(\frac{\hw_k}{\xdfw
_k}, \sum
_{0<\vert j\vert\leq k}\frac{\hw_j}{n\lambda_j} \biggr) \biggr]
\quad\mbox{and}
\nonumber
\\[-8pt]
\\[-8pt]
\psi_n& :=&\psi_n(\xdfw,\lambda,\hw):=\max\biggl(
\frac{\hw_{\kstar}}{\xdfw_{\kstar}}, \sum_{0<\vert j\vert\leq
\kstar}\frac{\hw_j}{n\lambda_j}
\biggr).
\nonumber
\end{eqnarray}
If, in addition,
$\eta:=\inf_{n\geq1}\{\psi_n^{-1}\min(\hw_{\kstar}\xdfw_{\kstar}^{-1},
\sum_{0<\vert l\vert\leq{\kstar}}{\hw_l}{(n\lambda_l})^{-1}) \}
>0$, then, for
all \mbox{$n\geq2$}
\[
\inf_{\tf} \sup_{\xdf\in\cF_\bw^\br} \bigl\{ \Ex\Vert\txdf-
\xdf
\Vert_\hw^2 \bigr\}\geq\frac{\eta\min(\xdfr-1,
1/(8\edfr\Gamma))}{16}
\psi_n,
\]
where the infimum is taken over all possible estimators of $f$.
\end{theo}
\begin{rem}
When $\edf$ is known, it is natural to consider the orthogonal series estimator
$\tf_k := 1 + \sum_{1< \vert j\vert\leq k} (\widehat{[\ydf]}_j /
[\edf
]_j) e_j$.
It is easily seen that for $\vert j\vert\leq k$, we have $\Ex[[\tf
]_j] =
[\xdf]_j$ and $\var([\tf]_j)\leq(n\vert[\edf]_j\vert^2)^{-1}$,
whereas $\Ex[[\tf]_j] = 0$
and $\var([\tf]_j) = 0$ for $\vert j\vert> k$. Thus, for all $\xdf
\in
\Fgr$ and
$\edf\in\Eld$, we have
\[
\Ex\bigl[\Vert\tf_k - \xdf\Vert_\hw^2
\bigr] \leq\sum_{\vert j\vert>k} \hw_j \bigl\vert[\xdf]_j\bigr\vert^2 +
\frac{1}{n}\sum
_{0< \vert j\vert\leq k} \frac{\hw_j}{\vert[\edf]_j\vert^2} \leq
(\xdfr+ \edfr) \max\biggl(
\frac{\hw_k}{\xdfw_k}, \sum_{0<\vert j\vert\leq k}\frac{\hw
_j}{n\lambda_j}
\biggr).
\]
Thus, the choice $\kstar$ of $k$ from \eqref{def:psin} realizes
the best variance--bias trade-off, $\psi_n$. This demonstrates that
when $\edf$ is known, $\tf_{\kstar}$
attains the rate $\psi_n$, which thus is minimax optimal.
\end{rem}
The proof of the last assertion is based on
Assouad's cube technique (Korostel{\"e}v and Tsybakov \cite
{KorostelevTsybakov1993}), which
involves constructing $2^{2\kstar}$ candidates of deconvolution
densities that have the largest possible $\Vert\cdot\Vert_\hw
$-distance but are still
statistically indistinguishable. Of note, the
additional assumption $\sum_{j\in\Z} \xdfw_j^{-1}=\Gamma<\infty$ is
used only to ensure that these candidates are densities. Also of note,
in the case where $\xdfr=1$, the lower bound is equal to 0,
because in this situation the set $\Fgr$ reduces to a singleton
containing only the uniform density. In the next theorem, we state a lower
bound characterizing the additional complexity due to the unknown
error density, which, surprisingly, depends only on the error sample
size.
\begin{theo}\label{theo:lower:m} Assume \eqref{eq:obs} and let $\hw
$, $\xdfw$, and $\lambda$ be sequences
satisfying Assumption~\textup{\ref{ass:minreg}}. For all
$m\geq2$, let
%
%
\begin{equation}
\label{def:kappam} \kappa_m:=\kappa_m(\xdfw,\lambda,
\hw):= \max_{j\in\N} \biggl\{ \hw_j\bw_j^{-1}
\min\biggl(1,\frac{1}{m\tw_j} \biggr)  \biggr\}.
\end{equation}
If in addition there exists a density in ${\mathcal E}^{\sqrt{\edfr
}}_\lambda$
that is bounded from below by $1/2$, then, for all $m\geq2$,
\[
\inf_{\tf}\sup_{\xdf\in\cF_\bw^\br}\sup_{\edf\in{\mathcal
E}_\tw
^\td} \bigl\{ \Ex\Vert
\txdf-\xdf\Vert_\hw^2 \bigr\}\geq\frac
{\min(\xdfr-1,1)\min(1/(4
d),(1-d^{-1/4})^2)}{4\sqrt d}
\kappa_m,
\]
where the infimum is taken over all possible estimators of $\xdf$.
\end{theo}
The proof of the last assertion takes its inspiration from a proof
given by Neumann \cite{Neumann1997}, who proved a similar lower bound for
deconvolution on the real line when both densities $\xdf$ and $\edf$
are ordinary smooth, that is, $\xdfw$ and $\lambda$ have polynomial
decay. In contrast to the proof of Theorem~\ref{theo:lower:n}, here we
only need compare two candidates of error
densities that are still statistically indistinguishable. However, to
ensure that these candidates are densities, we impose the additional
condition. It is easily seen that this condition is satisfied if
$\ell:=\sum_{j\in\Z} \lambda_j^{-1/2}<\infty$ and $\sqrt{d}\geq
\max(4\ell^2,1)$. Of note, in case where $\edfr=1$, the
set $\Eld$ of possible error densities reduces to a singleton, and
thus the lower bound is equal to 0. Finally, by a combination of
both lower bounds, we obtain the next corollary.

\begin{coro}\label{cor:lower}Under the assumptions of Theorem~\ref
{theo:lower:n} and~\ref{theo:lower:m} for all $n,m\geq2$
\[
\inf_{\tf}\sup_{\xdf\in\cF_\bw^\br}\sup_{\edf\in{\mathcal
E}_\tw
^\td} \bigl\{ \Ex\Vert
\txdf-\xdf\Vert_\hw^2 \bigr\}\geq C(\eta,\xdfr,\edfr,
\Gamma) \max(\psi_n,\kappa_m).
\]
\end{coro}

\subsection*{Upper bound} In the next theorem and all subsequent
results, we assume observations according to \eqref{eq:obs}.
First, we summarize sufficient conditions to ensure the optimality of
the orthogonal series estimator $\hxdf_{k}$ defined in \eqref{def:est},
provided that the dimension parameter $k$ is chosen appropriately. We use
the value $\kstar$ defined in \eqref{def:psin}, which, although
obviously involving the sequences $\hw, \xdfw$, and $\lambda$,
surprisingly does not depend on the $\epsilon$ sample size $m$. With
this choice, the estimator attains the lower bound given in
Corollary~\ref{cor:lower} up to a constant and thus is minimax-optimal.

\begin{theo}\label{theo:upper} Under Assumption~\textup{\ref
{ass:minreg}}, we
have, for all $n,m\geq1$,
\[
\sup_{\xdf\in\cF_\bw^\br}\sup_{\edf\in{\mathcal E}_\tw^\td}
\bigl\{ \Ex\Vert\hxdf_{\kstar}-
\xdf\Vert^2_\hw\bigr\}\leq C \bigl\{ (\edfr+ \xdfr)
\psi_n + \edfr\xdfr\kappa_{m} \bigr\}.
\]
\end{theo}
Note that under slightly stronger conditions on the sequences $\hw$,
$\xdfw$, and $\lambda$ than those in Assumption~\ref{ass:minreg}, it
can be shown
that in the case of equally large samples from $Y$ and $\epsilon$, we
always have the same rate as in the case of known error density.
However, in special cases, the required $\epsilon$ sample size can be
much smaller than the $Y$ sample size, as we show below.

\subsection*{Illustration: Estimation of derivatives}\label{sec:minimax:ill}
Here we illustrate our results considering classical smoothness
assumptions. Regarding the deconvolution density $\xdf$, it is
interesting to recall that the class $\Fgr$ is a
subset of the Sobolev space of $p$-times differentiable periodic
functions if
$\xdfw_j \sim\vert j\vert^{2p}$ (Neubauer \cite
{Neubauer88,Neubauer1988}). We call
this case
\textit{ordinary smooth}.
Moreover, up to a constant, for any function
$h\in\Fgr$, the weighted norm $\Vert h\Vert_\hw$ with $\hw_{j}\sim
j^{2s}$ equals the $L^2$ norm of the $s$th weak derivative $h^{(s)}$
for each integer $0\leq s\leq p$. By virtue of this relationship, the
results in the previous section imply both a lower bound and an upper
bound of the $L^2$ risk for estimation of the $s$th weak
derivative of $\xdf$.
If, in contrast, $\xdfw_j\sim
\exp(\vert j\vert^{2p})$ with $p>1$, then $\Fgr$ is a class of analytic
functions (Kawata \cite{Kwata1972}). We refer to this situation as
\textit{super smooth}.

As for the error densities, we consider two special cases
corresponding to a regular decay of their Fourier
coefficients. The error density is called \textit{ordinary smooth} if
$\lambda_j \sim\vert j\vert^{-2a}$ for some $a>1/2$ and \textit{super
smooth} if $\lambda_j\sim\exp(-\vert j\vert^{2a})$ for some $a>0$.

We consider the following three situations:
In the cases {[o-o]} and {[s-o]}, the error density is
ordinary smooth and the deconvolution density is either ordinary smooth
or super smooth case, respectively. Case {[o-s]} is the opposite
of case {[s-o]}.

It is readily seen that in all of these cases, the minimal regularity
conditions given in Assumption~\ref{ass:minreg} and the additional
conditions in
Theorems~\ref{theo:lower:n} and~\ref{theo:lower:m} translate to
simple restrictions on $p,a$, and $s$, which are given in the following
proposition. Roughly speaking, these restrictions imply that both the
deconvolution density and the error density are at
least continuous. The lower bounds presented in the following assertion
follow directly from Corollary~\ref{cor:lower}:
\begin{prop}\label{coro:ex:lower}
\begin{enumerate}[{[o-o]}]
\item[{[o-o]}] For $p>1/2$, $a>1$, and $0\leq
s \leq p$, we have for all $n,m\geq1$
\[
\inf_{\txdf^{(s)}}\sup_{\xdf\in\Fgr}\sup_{\edf\in{\mathcal
E}_\tw
^\td
} \bigl\{ \Ex\bigl\Vert
\txdf^{(s)}-\xdf^{(s)}\bigr\Vert^2 \bigr\}\gtrsim
n^{-2(p-s)/(2p+2a+1)} + m^{-((p-s)\wedge a)/a}.
\]
\item[{[s-o]}] For $p>0$, $a>1$, and $s\geq0$, we have for
all $n,m\geq1$
\[
\inf_{\txdf^{(s)}}\sup_{\xdf\in\Fgr}\sup_{\edf\in{\mathcal
E}_\tw
^\td
} \bigl\{ \Ex\bigl\Vert
\txdf^{(s)}-\xdf^{(s)}\bigr\Vert^2 \bigr\}\gtrsim
n^{-1} (\log n)^{(2a+2s+1)/(2p)} + m^{-1}.
\]
\item[{[o-s]}] For $p>1/2$, $a>0$,
and $0\leq s\leq p$, we have for all $n,m\geq1$
\[
\inf_{\txdf^{(s)}}\sup_{\xdf\in\Fgr}\sup_{\edf\in{\mathcal
E}_\tw
^\td
} \bigl\{ \Ex\bigl\Vert
\txdf^{(s)}-\xdf^{(s)}\bigr\Vert^2 \bigr\}\gtrsim(\log
n)^{-(p-s)/a} + (\log m)^{-(p-s)/a}.
\]
\end{enumerate}
\end{prop}
\begin{rem}
We do not treat the doubly exponential case [s-s] here,
because doing so would require rather intricate computations and
distinctions of
cases. A detailed analysis of this case in the context of density
deconvolution on the real line has been provided by Butucea and
Tsybakov \cite{BT:08a,BT:08b}.
Note that the expressions in $n$ in the foregoing result coincide with
the lower bounds for the deconvolution problem on the real line, which can
be found in the literature. For example, in cases where the error
distribution is known, Fan \cite{Fan91} have addressed the cases
{[o-o]} and \mbox{{[o-s]}} and Butucea
\cite{But:04} examined the case {[s-o]}.\footnote{When
comparing the bounds, attention must be given to the slightly
different parameterizations of the density classes in the cited
articles.} Those authors developed kernel-based
estimation procedures which attain these lower bounds. In the case
{[o-o]} (still on the real line), for cases where the error density
is unknown, \cite{Neumann1997} also investigated the impact of
estimating the error density and obtained the same lower bound as in
the foregoing result.
\end{rem}
As an estimator of $f^{(s)}$, we consider the $s$th weak derivative of
the estimator $\hxdf_k$ defined in
\eqref{def:est}, with $k$ as specified below. Given the exponential
basis $\{e_j\}_{j\in\Z}$, we recall
that for each integer $0\leq s\leq p$, the $s$th derivative in a weak
sense of the estimator $\hxdf_k$ is
%
%
\begin{equation}
\label{def:est:s} \hxdf^{(s)}_k=\sum
_{j\in\Z}(2\mathrm{i} \uppi j)^s[
\hxdf_k]_je_j.
\end{equation}
As an immediate consequence of
Theorem~\ref{theo:upper}, the rates of the lower bound given by
Proposition~\ref{coro:ex:lower} are attained for $k=\kstar$, as
summarized in the next result. Thus, we have proven that these rates
are optimal
and that the proposed estimator $\hxdf^{(s)}_{\kstar}$ is minimax
optimal in both
cases. Furthermore, it is of interest to characterize the minimal size
$m$ of the additional sample from $\epsilon$ needed to attain
the same rate as in case of a known error density. Thus, we let the
$\epsilon$-sample size depend on the
$Y$-sample size $n$ as well.

\begin{prop}\label{coro:ex:upper}Let $(m_n)_{n\geq1}$ be a sequence
of positive integers:
\begin{enumerate}[{[o-o]}]
\item[{[o-o]}] For $p>1/2$, $a>1$, and $0\leq
s \leq p$ with $\kstar\sim n^{1/(2p+2a+1)}$, we have for all $n,m\geq1$
\[
\sup_{\xdf\in\Fgr}\sup_{\edf\in{\mathcal E}_\tw^\td} \bigl\{
\Ex\bigl\Vert\hxdf^{(s)}_{\kstar}-
\xdf^{(s)}\bigr\Vert^2 \bigr\} \lesssim n^{-2(p-s)/(2p+2a+1)} +
m^{-( (p-s)\wedge a)/a}
\]
and if
$q_{\mathrm{o\mbox{-}o}}:=\lim_{n\to\infty}n^{2((p-s)\vee
a)/(2p+2a+1)} m_n^{-1}$ exists,\footnote{The limit ``$\infty$''
is authorized, with $\lim_{n\to\infty} a_n = \infty\dvtx\!\!\iff
\forall
K>0\ \exists
n_0\in\N\ \forall n\geq n_0 \dvt a_n \geq K$.} then it
follows that
as $n\to\infty$
\[
\sup_{\xdf\in\Fgr}\sup_{\edf\in
{\mathcal E}_\tw^\td} \bigl\{ \Ex\bigl\Vert\hxdf^{(s)}_{\kstar}-
\xdf^{(s)}\bigr\Vert^2 \bigr\}= %
\cases{ \mathrm{O}
\bigl(n^{-2(p-s)/(2p+2a+1)} \bigr) &\quad if $q_{\mathrm{o\mbox{-}o}} <
\infty$,
\cr
\mathrm{O} \bigl(m_n^{-((p-s)\wedge a)/a} \bigr)&\quad otherwise.}
\]
\item[{[s-o]}]For $p>0$, $a>1$, and $s\geq0$ with $\kstar
\sim
(\log n)^{1/(2p)}$, we have, for all $n,m\geq1$,
\[
\sup_{\xdf\in\Fgr}\sup_{\edf\in{\mathcal E}_\tw^\td} \bigl\{
\Ex\bigl\Vert\hxdf_{\kstar}^{(s)}-
\xdf^{(s)}\bigr\Vert^2 \bigr\} \lesssim n^{-1} (\log
n)^{(2a+2s+1)/(2p)} + m^{-1}
\]
and if $q_{\mathrm{s\mbox{-}o}}:=\lim_{n\to\infty} n (\log
n)^{-(2a+2s+1)/(2p)} m_n^{-1}$ exists, it
follows as
$n\to\infty$
\[
\sup_{\xdf\in\Fgr}\sup_{\edf\in
{\mathcal E}_\tw^\td} \bigl\{ \Ex\bigl\Vert\hxdf^{(s)}_{\kstar}-
\xdf^{(s)}\bigr\Vert^2 \bigr\}= %
\cases{ \mathrm{O}
\bigl(n^{-1} (\log n)^{(2a+2s+1)/(2p)} \bigr) &\quad if $q_{\mathrm{s\mbox{-}o}}
< \infty$,
\cr
\mathrm{O} \bigl(m_n^{-1} \bigr)&\quad
otherwise.} %
\]
\item[{[o-s]}] For $p>1/2$, $a>0$,
and $0\leq s\leq p$ with $\kstar\sim(\log n)^{1/(2a)}$, we
have, for all \mbox{$n,m\geq1$},
\[
\sup_{\xdf\in\Fgr}\sup_{\edf\in
{\mathcal E}_\tw^\td} \bigl\{ \Ex\bigl\Vert\hxdf^{(s)}_{\kstar}-
\xdf^{(s)}\bigr\Vert^2 \bigr\} \lesssim(\log n)^{-(p-s)/a} +
( \log m)^{-(p-s)/a}
\]
and if
$q_{\mathrm{o\mbox{-}s}} := \lim_{n\to\infty} (\log n) (\log
m_n)^{-1}$ exists, then it
follows, as $n\to\infty$,
\[
\sup_{\xdf\in\Fgr}\sup_{\edf\in
{\mathcal E}_\tw^\td} \bigl\{ \Ex\bigl\Vert\hxdf^{(s)}_{\kstar}-
\xdf^{(s)}\bigr\Vert^2 \bigr\}= %
\cases{\mathrm{O}
\bigl((\log n)^{-(p-s)/a} \bigr)&\quad if $q_{\mathrm{o\mbox{-}s}} < \infty$,
\cr
\mathrm{O} \bigl((\log m_n)^{-(p-s)/a} \bigr)&\quad otherwise.}
\]
\end{enumerate}
\end{prop}

The existence of the limits
$q_{\mathrm{o\mbox{-}o}}$, $q_{\mathrm{o\mbox{-}s}}$, and
$q_{\mathrm{s\mbox{-}o}}$ is required only to exclude the case of
oscillating sequences, which we are not interested in here. In this
case, none of the two terms in the upper
bound is asymptotically dominant, and the convergence rate is the
alternating maximum of the two terms.

In the case {[o-o]}, whenever $n^{2((p-s)\vee a)/(2p+2a+1)} =
\mathrm{O}(m_n)$, which is much less
than $m_n=n$, we obtain the rate of known error density. This is even
more visible in the case {[o-s]},
where the rate of known error density is attained even if $m_n=n^r$ for
arbitrarily small $r>0$. Moreover, we emphasize the influence of the
parameter $a$ that characterizes the rate of decay of the Fourier
coefficients of the error density $\edf$. Because a smaller value of $a$
leads to faster rates of convergence, this parameter is often called
\textit{degree of ill-posedness} (e.g., Natterer \cite{Natterer84}).

\section{Adaptive estimation}\label{sec:adaptive-known}

Our aim is to construct an adaptive
estimator of the deconvolution density $\xdf$. Adaptation means that
despite an unknown error density in $\Eld$, the estimator should
attain the optimal rate of convergence $\max(\psi_n, \kappa_m)$ over
the ellipsoid $\Fgr$ for a wide range of different weight sequences
$\xdfw$ and $\lambda$.

In a first step, we suppose that $\edf$ is known, but $\xdfw$ and
$\xdfr$ are unknown.
In what
follows, we consider the orthogonal series estimator $\hxdf_k$ defined in
\eqref{def:est} and construct a procedure to choose the dimension
parameter $k$ based on a model selection approach via penalization.
This partially adaptive choice $\tk$ will involve only the data and
the error density $\edf$.

In a second step, we replace $\edf$ with its empirical version and thus
dispense with any knowledge about $\edf$. Doing so, we obtain a fully
adaptive choice $\hk$ of the dimension parameter.

\subsection*{Partially adaptive estimation knowing {$\boldsymbol\edf$}}
\label{sec:patrially_adapdite}

We first introduce sequences that are used
below.

\begin{defin}\label{def:known}
For all $n,m\geq1$ and $k\geq0$, define
\begin{longlist}[(ii)]
\item
$ \Delta_k:=\Delta_k(\edf):={\max}_{-k\leq j\leq k}\frac{\hw
_j}{\vert[\edf]_j\vert^2}
\mbox{ and } \delta_k:=\delta_k(\edf) := 2 k
\Delta_k\frac{\log(\Delta_k\vee(k+2))}{\log
(k+2)}$;
\item given $\hw^+_k:=\max_{0\leq j\leq k}\hw_j$ and $\Cy_n^\circ:=
\max\{1\leq N\leq n \vert\hw^+_N\leq n\}$, let
\[
\Cy_n := \Cy_n(\edf) := \min\biggl\{1\leq j \leq
\Cy_n^\circ\Big\vert\frac{\vert[\edf]_j\vert^2}{j\hw^+_j} \leq
\frac
{\log(n+2)}{n}
\biggr\} - 1,
\]
defining further $b_m := (8\log( \log(m+20) )
)^{-1}$, let
\[
\Ce_m:= \Ce_m(\edf) := \min\bigl\{1\leq j \leq m \vert
\bigl\vert[\edf]_j\bigr\vert^2 \leq m^{-1 +b_m} \bigr\} - 1;
\]
\end{longlist}
with $\Cy_n := \Cy_n^\circ$ and $\Ce_m := m$ when the
respective set in the definition is empty.
\end{defin}
These sequences are used for small sample sizes as well, which
explains their rather complicated form.
We can now define a partially adaptive choice
of the dimension parameter $k$,
%
%
\begin{equation}
\label{gen:ada:1} \tk:= \mathop{\argmin}_{0\leq k\leq(\Cy_n\wedge
\Ce_m)} \biggl[-\Vert
\hxdf_{k}\Vert^2_\hw+ 60 \frac{\delta_k}{n}
\biggr],
\end{equation}
which obviously depends only on the data and the error
density $\edf$. We obtain the fully adaptive estimator
below by
introducing the empirical versions of $\delta, \Cy$, and $\Ce$ given in
Definition~\ref{def:known}.

For a fixed $\edf$, we could now derive an upper risk bound for the
partially adaptive estimator $\hxdf_{\tk}$, which would depend on
$\delta$, $\Cy$, and $\Ce$. But because we wish to obtain a uniform
upper risk bound over the class $\Eld$, instead we now
redefine the foregoing objects referring only to the weight sequence
$\lambda$ and the constant $\edfr$.

\begin{defin}\label{def:knownlambda} Let $\hw^+$, $N^\circ$, and $b$
as in Definition~\ref{def:known}.
\begin{longlist}[(iii)]
\item For all $k\geq0$, define
$\Delta^\lambda_k:=\max_{-k \leq j\leq k}\hw_j/\lambda_j$ and
\[
\label{deltamref} \delta_k^\lambda:= 2 k
\Delta^\lambda_k \frac{\log(\Delta^\lambda_k\vee(k+2))}{\log(k+2)}.
\]

\item Define two sequences, $\Cy^\lambda$ and
$\Ce^\lambda$, as follows:
\begin{eqnarray*}
\Cy^\lambda_n &:=& \min\biggl\{1\leq j \leq
\Cy_n^\circ\Big\vert\frac{\lambda_j}{j\hw_j^+} < \frac{4\edfr\log(n+2)}{n}
\biggr\} -1 ,
\\
\Ce^\lambda_m &:=& \min\bigl\{1\leq j \leq m \vert
\lambda_j < 4\edfr m^{-1 + b_m} \bigr\} -1.
\end{eqnarray*}
If the set in
the definition is empty, then we set $\Cy_n^\lambda:= 0$ or $\Ce
^\lambda_m := 0$,
respectively.

\item Define two sequences, $\Cy^u$ and $\Ce^u$, as follows:
\begin{eqnarray*}
\Cy_n^u := \Cy_n^u(\lambda) &:=&
\min\biggl\{1\leq j \leq n \Big\vert\frac{\lambda_j}{j\hw^+_j} <
\frac
{\log(n+2)}{4\edfr n} \biggr
\} - 1,
\\
\Ce_m^u:=\Ce_m^u(\lambda) &:=&
\min\biggl\{1\leq j \leq m \Big\vert\lambda_j < \frac{m^{-1 +
b_m}}{4\edfr}
\biggr\} -1.
\end{eqnarray*}
If the set in
the definition is empty, we set $\Cy_n^u := n$ or $\Ce^u_m := m$.

\item Let $\Sigma\dvtx\R\to\R$ be a non-decreasing function such that,
for all $C>0$,
\[
\label{ass:sum} \sum_{k\geq
1}C {\Delta_k^\lambda}
\exp\biggl(-\frac{k\log(\Delta_k^\lambda\vee(k+2))}{3 C \log
(k+2)} \biggr)\leq\Sigma(C)<\infty.
\]
\end{longlist}
\end{defin}

It is easy to see that there exists always a function $\Sigma$
satisfying the defining condition. Moreover, as we show in Lemma~\ref
{lem:schachtel} in the \hyperref[app]{Appendix}, the sequences defined
above satisfy
$\Cy^\lambda_n \leq\Cy_n\leq
\Cy^u_n$ and $\Ce^\lambda_m \leq\Ce_m \leq\Ce^u_m$ for all
$n,m\in
\N$. In the
illustration below we compute these objects explicitly.

\begin{theo}\label{theo:phi}
Let $\zeta_\edfr:= \log(3 \edfr)/\log(\edfr)$. Under
Assumption~\textup{\ref{ass:minreg}}, for all $n,m\geq1$,
\begin{eqnarray*}
\sup_{\xdf\in\cF_\xdfw^\xdfr} \sup_{\edf\in\Eld} \bigl\{ \Ex
\Vert\hxdf_{\tk}-
\xdf\Vert_\hw^2 \bigr\} &\leq& C \biggl\{ (\xdfr+ \edfr
\zeta_\edfr)\min_{0\leq k\leq(\Cy_n^\lambda\wedge
\Ce_m^\lambda)} \biggl[ \max\biggl( \frac{\hw_k}{\xdfw_{k}},
\frac{\delta_k^\lambda}{n} \biggr) \biggr] + \xdfr\edfr\kappa_m
\biggr\}
\\
&&{}+C(\xdfr,\edfr,\Edfw, \Sigma) \biggl[\frac{1}{m} + \frac{1}{n}
\biggr].
\end{eqnarray*}
\end{theo}

A comparison with the
lower bound from Corollary~\ref{cor:lower} shows that this upper bound
ensures minimax optimality of the estimator $\hxdf_{\tk}$ only if
\[
\label{eq:1}\psi_{n,m}^\diamond:=\min_{1\leq k\leq(\Cy_n^\lambda
\wedge\Ce
_m^\lambda)} \biggl[
\max\biggl( \frac{\hw_k}{\xdfw_{k}}, \frac{\delta_k^\lambda}{n}
\biggr) \biggr]
\]
is in the same order as $\psi_n
=\min_{k\in\N} \{\max(\frac{\hw_k}{\xdfw_k}, \sum_{0<\vert
j\vert\leq
k}\frac{\hw_j}{n\lambda_j} ) \}$. Note that, by construction,
$\delta_k^\lambda\geq\sum_{0<\vert j\vert\leq k}{\hw_j}{\lambda_j^{-1}}$
for all
$k\geq1$. In addition, $\delta^\lambda$ is direcly related to the penalty
function. The next assertion is a immediate consequence of
Theorem~\ref{theo:phi}, and we omit its proof.
\begin{coro}\label{coro:knownraten1}
Under Assumption~\textup{\ref{ass:minreg}}, and if
\[
\eta^\diamond:=\sup_{n,m\geq1} \bigl\{\psi_{n,m}^\diamond/
\max(\psi_n,\kappa_m) \bigr\}<\infty,
\]
then we have, for
all $n,m\geq1$,
\[
\sup_{\xdf\in\cF_\xdfw^\xdfr}\sup_{\edf\in{\mathcal E}_\lambda
^\edfr} \bigl\{ \Ex\Vert\hxdf_{\tk}-
\xdf\Vert_\hw^2 \bigr\} \leq C \bigl(\eta^\diamond,
\Sigma,\xdfr,\edfr,\Edfw\bigr) \max(\psi_n, \kappa_m).
\]
\end{coro}

In Theorem~\ref{theo:upper}, we have shown the minimax optimality of
the orthogonal series estimator under the optimal choice $\kstar$ of the
dimension parameter. Comparing Corollary~\ref{coro:knownraten1} with
this theorem, it is noteworthy that the only additional assumption
needed to ensure minimax optimality of the partially adaptive
estimator is $\eta^\diamond<\infty$.

\begin{rem}\label{rem:klambda}
The partially adaptive choice $\tk$ still
depends on $\edf\in\Eld$.
However, we can already
define a procedure depending only on the sequence $\lambda$ and
the constant $\edfr$, namely
\[
\label{gen:ada:1edfw} \tk^\lambda:= \mathop{\argmin}_{1\leq k\leq
(\Cy
_n^\lambda\wedge
\Ce_m^\lambda)} \biggl[-
\Vert\hxdf_{k}\Vert^2_\hw+ 60
\frac{d \delta^\lambda_k}{n} \biggr].
\]
Roughly speaking, this choice requires knowledge of the degree of
ill-posedness of the underlying inverse problem only. It is
straightforward to derive an upper risk bound for $\hxdf_{\tk^\lambda}$,
which is, up to minor changes in the constants, the same as that in
Theorem~\ref{theo:phi}. Its proof follows the lines of the proof of
Theorem~\ref{theo:phi}, using the new penalty term $\pen(k) =
60 \edfr\delta_k^\lambda$. The only change occurs when applying
Lemma~\ref{lem:talalem}, which uses
$\delta_k^\ast=\edfr\delta_k^\lambda$ and
$\Delta_k^\ast=\edfr\Delta_k^\lambda$ rather than
$\delta_k^\ast=\delta_k$ and $\Delta_k^\ast=\Delta_k$.
\end{rem}

\subsection*{Fully adaptive estimation}\label{sec:adaptive:unknown}

We begin by defining empirical versions of the sequences given
in Definition~\ref{def:known}.

\begin{defin}\label{def:unknown} For all $n,m\geq1$ and $k\geq0$, define
\begin{longlist}[(ii)]
\item
$\hDelta_k := {\max}_{-k\leq j \leq k}\frac{\hw_j}{\vert\widehat
{[\edf ]}_j\vert^{2}}\mathbh{1}{\{\vert\widehat{[\edf]}_j\vert
^2\geq1/m\}}$ and
$\widehat{\delta}_k := k \hDelta_k \frac{\log(\hDelta_k \vee
(k+2))}{\log(k+2)}$;
\item given $\Cy_n^\circ$, $\hw^+$, and $b$ from Definition~\ref
{def:known},
\begin{eqnarray*}
\hCy_n &:=& \min\biggl\{1\leq j\leq\Cy_n^\circ
\Big\vert\frac{\min(\vert\widehat{[\edf]}_j\vert^2,\vert\widehat
{[\edf ]}_{-j}\vert^2)}{j\hw^+_j} < \frac{\log
(n+2)}{n} \biggr\} - 1,
\\
\hCe_m &:=& \min\bigl\{1 \leq j\leq m \vert\min\bigl(\bigl\vert
\widehat{[\edf]}_j\bigr\vert^2,\bigl\vert\widehat{[\edf]}_{-j}
\bigr\vert^2 \bigr) < m^{-1+b_m} \bigr\} - 1,
\end{eqnarray*}
with $\hCy_n := \Cy^\circ_n$ and $\hCe_m := m$ if the
respective sets in the definition are empty.
\end{longlist}
\end{defin}
We now define a data-driven choice of $k$, which, in contrast to $\tk
$, depends not on the sequences $\delta$, $\Cy$, or $\Ce$, but
rather on
$\widehat{\delta}$, $\hCy$, and $\hCe$:
%
%
\begin{equation}
\label{eq:8} \hk:= \mathop{\argmin}_{0\leq k\leq(\hCy_n\wedge
\hCe_m)} \biggl[ -\Vert\hf_k
\Vert^2_\hw+ 600 \frac{\widehat{\delta}_k}{n} \biggr].
\end{equation}
The constant 600 arising in the definition of
$\hk$, although convenient for deriving the theory, may be far too large
in practice and instead be determined by means of a simulation study,
as done by Comte, Rozenholc and Taupin \cite
{ComteRozenholcTaupin2006}, for example.

In the proof of Theorem~\ref{theo:phi}, we used\vspace*{1pt}
$(\Cy^\lambda_n\wedge\Ce^\lambda_m) \leq(\Cy_n\wedge\Ce_m)\leq
(\Cy^u_n\wedge\Ce^u_m)$ (Lemma~\ref{lem:schachtel}). In the proof of
the next theorem, we consider the event $\{(\Cy^\lambda_n\wedge
\Ce^\lambda_m) \leq(\hCy_n\wedge\hCe_m)\leq(\Cy^u_n\wedge\Ce
^u_m)\}$, on
which we can imitate the proof of Theorem~\ref{theo:phi}. To
control the risk on the complement of this event, we need to bound its
probability, which necessitates the following assumption.

\begin{assumption}\label{ass:MplusEins}
Suppose that
$ m^7\exp( - {m \lambda_{\Ce^u_m+1}}/({72 \edfr})
) \leq C(\lambda,\edfr)$ for all $m\geq1$.
\end{assumption}

\begin{theo}\label{thm:upper:unknown}
Under Assumptions~\textup{\ref{ass:minreg}} and~\textup{\ref
{ass:MplusEins}}, we have,
for all $n,m\geq1$,
\begin{eqnarray*}
\sup_{\xdf\in\cF_\xdfw^\xdfr} \sup_{\edf\in\Eld} \bigl\{ \Ex
\Vert\hxdf_{\hk}-
\xdf\Vert_\hw^2 \bigr\} &\leq& C \biggl\{ (\xdfr+ \edfr
\zeta_\edfr)\min_{0\leq k\leq(\Cy_n^\lambda\wedge
\Ce_m^\lambda)} \biggl[ \max\biggl( \frac{\hw_k}{\xdfw_{k}},
\frac{\delta_k^\lambda}{n} \biggr) \biggr] + \xdfr\edfr\kappa_m
\biggr\}
\\
&&{}+C(\xdfr,\edfr,\lambda,\Sigma) \biggl[\frac{1}{m} + \frac{1}{n}
\biggr].
\end{eqnarray*}
\end{theo}

\begin{rem}\label{rem:full_adapt}
Up to a change in the constant in front of the negligible terms, we
obtain the same bound as for the partially adaptive estimator
(Theorem~\ref{theo:phi}). Compared with Theorem~\ref{theo:phi}, the
only additional
assumption is \ref{ass:MplusEins}.
Note that in Lemma~\ref{lem:NuMu}(ii) in the \hyperref
[app]{Appendix}, we show that
$ m^7\exp( - {m \lambda_{\Ce^u_m}}/({72 \edfr})
) \leq C(\edfr)$ for all $m\geq1$ using only Assumption~\ref{ass:minreg}.
However, it is not obvious to us that Assumption~\ref{ass:minreg} also
implies the slightly stronger assertion $ m^7\exp( - {m \lambda_{\Ce
^u_m+1}}/({72 \edfr})
) \leq C(\edfr)$ for sufficiently large $m$, although in the
illustrations below, we show that Assumption~\ref{ass:MplusEins} is satisfied.
\end{rem}
Comparing Theorem~\ref{thm:upper:unknown} with the lower bound from
Corollary~\ref{cor:lower} shows
that this upper bound does not necessarily ensure minimax optimality of
the estimator
$\hxdf_{\hk}$. However, as in the partially adaptive case (cf.
Corollary~\ref{coro:knownraten1}), under the
additional assumption $\eta^\diamond<\infty$, the next assertion
establishes its optimality. Because this is an immediate consequence of
Theorem~\ref{thm:upper:unknown}, we omit the proof.

\begin{coro}\label{coro:knownraten} Under Assumptions~\textup{\ref
{ass:minreg}}
and~\textup{\ref{ass:MplusEins}}, and if
\[
\eta^\diamond=\sup_{n,m\geq1} \bigl\{\psi_{n,m}^\diamond/
\max(\psi_n,\kappa_m) \bigr\}<\infty,
\]
we have, for all $n,m \geq1$,
\[
\sup_{\xdf\in\cF_\xdfw^\xdfr}\sup_{\edf\in{\mathcal E}_\lambda
^\edfr} \bigl\{ \Ex\Vert\hxdf_{\hk}-
\xdf\Vert_\hw^2 \bigr\} \leq C \bigl(\eta^\diamond,
\Sigma,\xdfr,\edfr,\Edfw\bigr) \max(\psi_n, \kappa_m) .
\]
\end{coro}

\subsubsection*{Conclusion} The minimax optimality of the estimator
$\hxdf_{\kstar}$ has been shown under Assumption~\ref{ass:minreg} in
Theorem~\ref{theo:upper}, where the choice $\kstar$ of the dimension
parameter depends on the deconvolution density $\xdf$ and the error
density $\edf$. We have developed a fully data-driven
choice $\hk$. The foregoing results show that we need only the
additional Assumptions~\ref{ass:MplusEins} and $\eta^\diamond<\infty$
for the adaptive estimator $\hxdf_{\hk}$ to be minimax
optimal as well.

\subsection*{Illustration: Estimation of derivatives (continued from
Section~\protect\ref{sec:minimax})}
The following result shows that
without any prior knowledge on the error density $\edf$, the
adaptive penalized estimator automatically attains the optimal
rate in the cases {[o-s]} and {[s-o]} and in the case
{[o-o]} if $p-s> a$. Recall that the computation of the dimension
parameter $\hk$ given in \eqref{eq:8} involves the sequence
$\Cy^\circ$, which in our illustration satisfies $\Cy_n^\circ
\asymp n^{1/(2s)}$.

\begin{prop}\label{prop:ex-cont-ada-ukn}
Let $(m_n)_{n\geq1}$ be a sequence of positive integers and suppose
that the limits $q_{\mathrm{o\mbox{-}o}}$,
$q_{\mathrm{o\mbox{-}s}}$, and $q_{\mathrm{s\mbox{-}o}}$ defined in
Proposition~\ref{coro:ex:upper} exist in the respective cases.
\begin{enumerate}[{[o-o]}]
\item[{[o-o]}] We have that
\begin{eqnarray*}
\Delta_k^\lambda&\sim& k^{2a+2s},\qquad
\delta_k^\lambda\asymp k^{2a+2s+1},\qquad
\psi_{n,{m_n}}^\diamond\sim\bigl(\kstar\wedge\Ce_{m_n}^\lambda
\bigr)^{-2(p-s)},
\\
\Cy_n^\lambda&\asymp&(n/\log n)^{1/(2a+2s+1)},\qquad
\Ce_{m_n}^\lambda\asymp m_n^{(1-b_m)/(2a)}.
\end{eqnarray*}
In the case wehere $p-s>a$, the adaptive estimator $\hxdf^{(s)}_{\hk
}$ attains the
optimal rates (see Proposition~\ref{coro:ex:upper}).
In the case where $p-s\leq a$, if $q_{\mathrm{o\mbox{-}o}}<\infty$,
then we have,
supposing that $q^b_{\mathrm{o\mbox{-}o}}:=\lim_{n\to\infty}
n^{2a/(2p+2a+1)} m_n^{-1+b_{m_n}}$ exists,
\[
\sup_{\xdf\in\Fgr}{\sup_{\edf\in\Eld}} \bigl\{ \Ex\bigl\Vert\hxdf
^{(s)}_{\hk}-
\xdf^{(s)}\bigr\Vert^2 \bigr\}= %
\cases{ \mathrm{O}
\bigl(n^{-2(p-s)/(2p+2a+1)} \bigr) &\quad if $q^b_{\mathrm{o\mbox{-}o}} <
\infty$,
\cr
\mathrm{O} \bigl(m_n^{-(p-s)/a} m_n^{b_{m_n}}
\bigr)&\quad otherwise,} %
\]
whereas if $q_{\mathrm{o\mbox{-}o}} = \infty$, then we have
\[
\sup_{\xdf\in\Fgr}{\sup_{\edf\in\Eld}} \bigl\{ \Ex\bigl\Vert\hxdf
^{(s)}_{\hk}-
\xdf^{(s)}\bigr\Vert^2 \bigr\}= \mathrm{O} \bigl(m_n^{-(p-s)/a}m_n^{b_m}
\bigr).
\]
\item[{[s-o]}] The sequences $\Delta^\lambda$, $\delta^\lambda$,
$\Cy^\lambda$, and $\Ce^\lambda$ are the same as above. We have that
$\psi^\diamond_{n,m_n} \sim
(\kstar\wedge\Ce_{m_n}^\lambda)^{2s}\exp(-(\kstar\wedge\Ce
_{m_n}^\lambda)^{2p})$,
and $\hxdf^{(s)}_{\hk}$ attains the optimal rates.

\item[{[o-s]}] We have that
\begin{eqnarray*}
\Delta_k^\lambda&=& k^{2s}\exp
\bigl(k^{2a} \bigr),\qquad\delta_k^\lambda\asymp
k^{2a+2s+1}\exp\bigl(k^{2a} \bigr) (\log k)^{-1},
\\
\psi_{n,m}^\diamond&\sim& \bigl(\kstar\wedge
\Ce_{m_n}^\lambda\bigr)^{-2(p-s)},
\\
\Cy_n^\lambda&\asymp& \bigl(\log\bigl(n/(\log
n)^{(2a+2s+1)/(2a)} \bigr) \bigr)^{1/(2a)},\qquad\Ce_{m_n}^\lambda
\asymp\bigl(({1-b_m})\log m_n \bigr)^{1/(2a)},
\end{eqnarray*}
and the adaptive estimator $\hxdf^{(s)}_{\hk}$ attains the optimal rates.
\end{enumerate}
\end{prop}

The adaptive estimator always attains the minimal rates
if $n \lesssim m_n$. We emphasize that this still holds when
$m_n\lesssim n$, except in the case {[o-s]} when the error
density is smoother than the $s$th derivative of the deconvolution
density ($p-s\leq a$) and when at the same time $m_n$ grows far more
slowly than $n$. The estimation of $\edf$ is negligible as soon as
$m_n^{1-b_{m_n}}$
grows at least as fast as $n^{2a/(2p+2a+1)}$ in this situation, whereas
in the nonadaptive case, only $m_n$ must satisfy this condition.
In the lossy case, the convergence rate
differs from the optimal rate by a factor $m_n^{b_{m_n}}$ only;
however, the
exponent $b_{m_n}$ tends to 0 as $n$ tends to infinity.

If considering the {[o-s]} case only, we could replace the bound
$m^{-1+b_m}$ by $m^{-1}\log m$ in the definition of
$\Ce^u$ (Definition~\ref{def:knownlambda}). Using this definition,
Assumption~\ref{ass:MplusEins} would still hold, and applying
Theorem~\ref{thm:upper:unknown}, the adaptive estimator would miss the
optimal rates by a logarithmic factor in the lossy case only.
However,
Assumption~\ref{ass:MplusEins} is violated in the super smooth case
under this definition of $\Ce^u$.

\begin{appendix}\label{app}
\section*{Appendix: Proofs}
\subsection{\texorpdfstring{Proofs of Section~\protect\ref{sec:minimax} (minimax theory)}
{Proofs of Section 2 (minimax theory)}}\label{app:proofs:minimax}
\subsubsection*{Lower bounds}

\begin{pf*}{Proof of Theorem~\ref{theo:lower:n}} Given $\zeta:=\eta
\min(\xdfr-1,
1/(8\edfr\Gamma))$ and $\alpha_n:=\psi_n (\sum_{0<\vert j\vert
\leq\kstar}
\hw_j/\break(\lambda_j\* n))^{-1}$, we consider the function $\xdf:= 1+
(\zeta\alpha_n/n)^{1/2}\sum_{0<\vert j\vert\leq\kstar}\lambda_j^{-1/2}
\ef_j$. We show that for any
$\theta:=(\theta_j)\in\{-1,1\}^{2\kstar}$, the function $\xdf
_\theta:=
1+ \sum_{0<\vert j\vert\leq\kstar} \theta_j[\xdf]_j \ef_j$
belongs to\vspace*{1pt}
$\cF_\xdfw^\xdfr$ and thus is a possible candidate for the
deconvolution density. For each $\theta$, the $Y$ density
corresponding to the $X$ density $\xdf_\theta$ is given by
$\ydf_\theta:=\xdf_\theta\ast\edf$. We denote by $\ydf_\theta
^n$ the
joint density of an i.i.d. $n$ sample from $\ydf_\theta$ and by
$\Ex_\theta$ the expectation with respect to the joint density
$\ydf_\theta^n$. Furthermore, for $0<\vert j\vert\leq\kstar$ and each
$\theta$,
we introduce $\theta^{(j)}$ by $\theta^{(j)}_{l}=\theta_{l}$ for
$j\ne
l$ and $\theta^{(j)}_{j}=-\theta_{j}$. The key argument of this proof
is the following reduction scheme. If $\txdf$ denotes an estimator of
$\xdf$, then we conclude that
\begin{eqnarray*}
\sup_{\xdf\in\Fgr} \Ex\Vert\txdf-\xdf\Vert_\hw^2 &
\geq& \sup_{\theta\in
\{-1,1\}^{2\kstar}} \Ex_\theta\Vert\txdf-\xdf_\theta
\Vert_\hw^2 \geq\frac{1}{2^{2{\kstar}}}\sum
_{\theta\in\{-1,1\}^{2\kstar}}\Ex_\theta\Vert\txdf-\xdf_\theta
\Vert_\hw^2
\\
&\geq&\frac{1}{2^{2{\kstar}}}\sum_{\theta\in\{-1,1\}^{2\kstar
}}\sum
_{0<\vert j\vert\leq\kstar}\hw_j\Ex_{{\theta}}\bigl\vert[\txdf-
\xdf_\theta]_j\bigr\vert^2
\\
&=& \frac{1}{2^{2{\kstar}}}\sum_{0<\vert j\vert\leq\kstar}\frac
{\hw_j}{2}
\sum_{\theta\in
\{-1,1\}^{2\kstar}} \bigl\{\Ex_{{\theta}}\bigl\vert[\txdf-
\xdf_\theta]_j\bigr\vert^2+\Ex_{{\theta^{(j)}}}\bigl\vert[
\txdf- \xdf_{\theta^{(j)}}]_j\bigr\vert^2 \bigr\},
\end{eqnarray*}
where for each $0<\vert j\vert\leq\kstar$ and any
function $F\dvtx\{-1,1\}^{2\kstar} \to\R$, we have
\[
\sum_{\theta\in\{-1,1\}^{2\kstar}}F(\theta) = \sum
_{\theta\in\{-1,1\}^{2\kstar}}F \bigl(\theta^{(j)} \bigr).
\]
Below we show that for all $n\geq2$, we have
%
\setcounter{equation}{0}
\renewcommand{\theequation}{A.\arabic{equation}}
\begin{equation}
\label{eq:443} \bigl\{\Ex_{{\theta}}\bigl\vert[\txdf-\xdf_\theta]_j
\bigr\vert^2+ \Ex_{{\theta^{(j)}}}\bigl\vert[\txdf-\xdf_{\theta^{(j)}}]_j
\bigr\vert^2 \bigr\} \geq\frac{\zeta\alpha_n}{4\lambda_j n}.
\end{equation}
Combining the last lower bound and the reduction scheme gives
\[
\sup_{\xdf\in\Fgr} \Ex\Vert\txdf-\xdf\Vert^2_\hw\geq
\frac{1}{2^{2{\kstar}}}\sum_{\theta\in\{-1,1\}^{2\kstar}}\sum
_{0<\vert j\vert\leq\kstar}\frac{\hw_j}{2}\frac{\zeta\alpha
_n}{4\lambda_j
n} =
\frac{\zeta}{8} \alpha_n\sum_{0<\vert j\vert\leq\kstar}
\frac{\hw_j}{\lambda_j n}.
\]
Thus, using the definition of $\zeta$ and $\alpha_n$, we obtain
the lower bound given in the theorem.

To conclude the proof, it remains to check \eqref{eq:443} and
$\xdf_\theta\in\Fgr$ for all $\theta\in\{-1,1\}^{2\kstar}$. The
latter is easily verified if $\xdf\in\Fgr$. To show that
$\xdf\in\Fgr$, we first note that $\xdf$ integrates to
1. Moreover, $f$ is nonnegative, because $\vert\sum_{0<\vert j\vert
\leq\kstar}
[f]_j\ef_j\vert\leq1$ and $\Vert\xdf\Vert_\xdfw^2\leq\xdfr$, which
can be
realized as follows. Using the condition
{$\sum_{j\in\Z}\xdfw_j^{-1}=\Gamma<\infty$}, we have
\begin{eqnarray*}
\biggl\vert\sum_{0<\vert j\vert\leq\kstar} [f]_j
\ef_j\biggr\vert&\leq&\sum_{0<\vert j\vert\leq\kstar} \bigl\vert
[f]_j\bigr\vert= \biggl( \frac{\zeta\alpha_n}{n} \biggr)^{1/2}\sum
_{0<\vert j\vert\leq\kstar} \lambda_j^{-1/2}
\\
&\leq&(\zeta\alpha_n )^{1/2} \biggl(\sum
_{0<\vert j\vert\leq\kstar} \xdfw_j^{-1}
\biggr)^{1/2} \biggl(\sum_{0<\vert j\vert\leq\kstar}
\frac{\xdfw_j}{n\lambda_j} \biggr)^{1/2}
\\
&\leq&(\zeta\alpha_n\Gamma)^{1/2} \biggl(\sum
_{0<\vert j\vert\leq\kstar} \frac{\xdfw_j}{n\lambda_j} \biggr)^{1/2}.
\end{eqnarray*}
Because
$\hw/\xdfw$ is nonincreasing, the definitions of $\zeta$, $\alpha_n$,
and $\eta$ imply that
%
%
\begin{equation}
\label{pr:theo:lower:n:e1} \biggl\vert\sum_{0<\vert j\vert\leq\kstar}
[f]_j\ef_j\biggr\vert\leq(\zeta\Gamma)^{1/2}
\biggl( \frac{\xdfw_{\kstar}}{\hw_{\kstar}} \alpha_n\sum
_{0<\vert j\vert\leq
\kstar}
\frac{\hw_j}{\lambda_j n} \biggr)^{1/2}\leq{ \biggl(\frac{\zeta
\Gamma}{\eta}
\biggr)^{1/2}\leq1},
\end{equation}
as well as $\Vert\xdf\Vert_\xdfw^2 \leq
1+\zeta\frac{\xdfw_{\kstar}}{\hw_{\kstar}}\alpha_n (\sum
_{0<\vert j\vert\leq\kstar}
\frac{\hw_j}{n\lambda_j} ) \leq{ 1+\zeta/\eta\leq r}$.

It remains to show \eqref{eq:443}. Consider the Hellinger affinity
$\rho(\ydf_{\theta}^n,\ydf_{\theta^{(j)}}^n)= \int
\sqrt{\ydf_{\theta}^n}\sqrt{\ydf_{\theta^{(j)}}^n}$. We then
obtain that,
for any estimator $\txdf$ of $\xdf$,
\begin{eqnarray*}
\rho\bigl(\ydf_{\theta}^n,\ydf_{\theta^{(j)}}^n
\bigr) &\leq&\int\frac{\vert[\txdf-\xdf_{\theta
^{(j)}}]_j\vert}{\vert[\xdf_\theta-\xdf_{\theta^{(j)}}]_j\vert
}\sqrt{\ydf_{\theta^{(j)}}^n}
\sqrt{\ydf_{\theta}^n} + \int\frac{\vert[\txdf-\xdf_{\theta
}]_j\vert}{\vert[\xdf_\theta-\xdf_{\theta^{(j)}}]_j\vert} \sqrt
{\ydf_{\theta}^n}\sqrt{\ydf_{\theta^{(j)}}^n}
\\
&\leq& \biggl( \int\frac{\vert[\txdf-\xdf_{\theta
^{(j)}}]_j\vert^2}{\vert[\xdf_\theta-\xdf_{\theta^{(j)}}]_j\vert
^2} \ydf_{\theta^{(j)}}^n
\biggr)^{1/2} + \biggl( \int\frac{\vert[\txdf-\xdf_{\theta
}]_j\vert^2}{\vert[\xdf_\theta-\xdf_{\theta^{(j)}}]_j\vert^2}
\ydf_{\theta}^n
\biggr)^{1/2}.
\end{eqnarray*}
Rewriting the last estimate, we obtain
%
%
\begin{equation}
\label{pr:theo:lower:n:e2} \bigl\{\Ex_{{\theta}}\bigl\vert[\txdf-
\xdf_\theta]_j\bigr\vert^2+\Ex_{{\theta^{(j)}}}\bigl\vert[
\txdf- \xdf_{\theta^{(j)}}]_j\bigr\vert^2 \bigr\} \geq
\tfrac{1}{2} \bigl\vert[\xdf_\theta-\xdf_{\theta^{(j)}}]_j
\bigr\vert^2 \rho^2 \bigl(\ydf_{\theta}^n,
\ydf_{\theta^{(j)}}^n \bigr).
\end{equation}
Next, we bound from below the Hellinger affinity $\rho(\ydf_{\theta
}^n,\ydf_{\theta^{(j)}}^n)$. Therefore, we first consider the
Hellinger distance,
\begin{eqnarray*}
H^2(\ydf_{\theta},\ydf_{\theta^{(j)}}) &:=& \int(\sqrt
\ydf_{\theta}-\sqrt\ydf_{\theta^{(j)}} )^2
\\
&\hspace*{2.8pt}=& \int\frac{ \vert\ydf_{\theta}-\ydf_{\theta
^{(j)}} \vert^2}{ (\sqrt\ydf_{\theta}+\sqrt\ydf_{\theta^{(j)}}
)^2}\leq4 \Vert\ydf_{\theta
}-
\ydf_{\theta^{(j)}} \Vert^2 = 16 \bigl\vert[\xdf]_j
\bigr\vert^2\bigl\vert[ \edf]_j\bigr\vert^2 \leq{
\frac{16\zeta\edfr}{\eta n}},
\end{eqnarray*}
where we have used that $\alpha_n\leq1/\eta$, $\edf\in
{\mathcal E}_\lambda^\edfr$, and $\ydf_\theta\geq1/2$ because
$\vert\sum_{0<\vert j\vert\leq\kstar} [\ydf_\theta]_j\ef_j\vert
\leq1/2$, which
can be
realized as follows. Using the condition {
$\sum_{j\in\Z}\xdfw_j^{-1}=\Gamma<\infty$} and $\edf\in
{\mathcal E}_\lambda^\edfr$, we obtain, in analogy to the proof of
\eqref{pr:theo:lower:n:e1}, that
\[
\biggl\vert\sum_{0<\vert j\vert\leq\kstar} [\ydf_\theta]_j
\ef_j\biggr\vert\leq\sum_{0<\vert j\vert\leq\kstar} \bigl\vert[
\xdf]_j\bigr\vert\bigl\vert[\edf]_j\bigr\vert\leq\biggl(
\frac{\zeta\alpha_n\edfr}{n} \biggr)^{1/2} \sum_{0<\vert j\vert
\leq\kstar}
\lambda_j^{-1/2}\leq{ \biggl(\frac{\zeta\edfr\Gamma}{\eta}
\biggr)^{1/2}\leq1/2}.
\]

Therefore, the definition of $\zeta$ implies
$H^2(\ydf_{\theta},\ydf_{\theta^{(j)}}) \leq2/n$. Using the
independence, that is,
$\rho(\ydf_{\theta}^n,\ydf_{\theta^{(j)}}^n)=\rho(\ydf_{\theta
},\ydf_{\theta^{(j)}})^n$,
together with the identity
$\rho(\ydf_{\theta},\ydf_{\theta^{(j)}})=1-\frac{1}{2}H^2(\ydf
_{\theta},\ydf_{\theta^{(j)}})$,
it follows $\rho(\ydf_{\theta}^n,\ydf_{\theta^{(j)}}^n) \geq(1-
n^{-1})^n\geq1/4$ for all $n\geq2$. Combining the last
estimate with \eqref{pr:theo:lower:n:e2}, we obtain \eqref{eq:443},
which completes the proof.
\end{pf*}

\begin{pf*}{Proof of Theorem~\ref{theo:lower:m}} We construct for each
$\theta\in\{-1,1\}$ an error density $\edf_\theta\in{\mathcal
E}_\lambda^\edfr$
and a deconvolution density $\xdf_\theta\in\cF_\xdfw^\xdfr$, such that
$\ydf_\theta:=\xdf_\theta\ast\edf_\theta$ satisfies $\ydf_1 =
\ydf_{-1}$. To be more precise, define $\kstarm:=\argmax_{\vert
j\vert>0}
[\hw_j\xdfw_j^{-1} \min(1,m^{-1}\lambda_j^{-1})]$ and $\alpha
_m:=\zeta
\min(1,m^{-1/2}\lambda_{\kstarm}^{-1/2})$ with $\zeta:=\min
(1/(2\sqrt
d),(1-d^{-1/4}))$. Observe that $1\geq(1-\alpha_m)^2\geq
(1-(1-1/d^{1/4}))^2\geq1/d^{1/2}$ and $1\leq(1+\alpha_m)^2\leq
(1+(1-1/d^{1/4}))^2=(2-1/d^{1/4})^2 \leq d^{1/2}$, which implies
$1/d^{1/2} \leq(1+\theta\alpha_m)^2\leq d^{1/2}$. We use these inequalities
below without further reference. By assumption, there is a
density $\edf\in{\mathcal E}_\lambda^{\sqrt\edfr}$ such that $\edf
\geq
1/2$. We
show below that for each $\theta$, the function
$\xdf_\theta:=1+(1-\theta\alpha_m) \frac{\min(\sqrt{r-1}
,1)}{d^{1/4}} \xdfw_{\kstarm}^{-1/2} e_{\kstarm}$ belongs to $\Fgr$,
and that the function $\edf_\theta:=\edf+ \theta\alpha_m [\edf
]_{\kstarm}
e_{\kstarm}$ is an element of ${\mathcal E}_\lambda^{\sqrt\edfr}$.
Moreover, it
is easily verified that $ \ydf_\theta= 1+ (1-\alpha_m^2)
\frac{\min(\sqrt{r-1},1)}{d^{1/4}} \xdfw_{\kstarm}^{-1/2}
[\edf]_{\kstarm} e_{\kstarm}$, and thus $g_1=g_{-1}$. We denote by
$\ydf_\theta^n$ the joint density of an i.i.d. $n$-sample from
$\ydf_\theta$ and by $\edf_\theta^m$ the joint density of an
i.i.d.
$m$ sample from $\edf_\theta$. Because the samples are independent of
one another, $p_\theta:=\ydf_\theta^n\edf_\theta^m$ is the joint
density of all observations, and we denote by $\Ex_\theta$ the
expectation with respect to $p_\theta$. Applying a reduction scheme,
we deduce
that for each estimator $\txdf$ of~$\xdf$,
\[
\sup_{\xdf\in\cF_\xdfw^\xdfr}\sup_{\edf\in{\mathcal E}_\lambda
^\edfr} \Ex\Vert\txdf-\xdf
\Vert^2_\hw\geq\max_{\theta\in\{-1,1\}} \Ex_{{\theta}}
\Vert\txdf-\xdf_\theta\Vert^2_\hw\geq
\frac{1}{2} \bigl\{\Ex_{1}\Vert\txdf- \xdf_1
\Vert^2_\hw+\Ex_{-1}\Vert\txdf-
\xdf_{-1}\Vert^2_\hw\bigr\}.
\]
Below we also show that for all $m\geq2$, we have
%
%
\begin{equation}
\label{pr:theo:lower:m:e1} \Ex_{1}\Vert\txdf-\xdf_1
\Vert^2_\hw+ \Ex_{-1}\Vert\txdf-
\xdf_{-1}\Vert^2_\hw\geq\tfrac{1}{8} \Vert
\xdf_{1}-\xdf_{-1}\Vert^2_\hw.
\end{equation}
Moreover, we have $ \Vert\xdf_1-\xdf_{-1}\Vert^2 = 4\alpha_m^2
\hw_{\kstarm}\xdfw_{\kstarm}^{-1} \frac{(\xdfr-1)\wedge
1}{d^{1/2}} =
4\frac{(\xdfr-1)\wedge1}{d^{1/2}} \zeta^2
\hw_{\kstarm}\xdfw_{\kstarm}^{-1}
\min(1,\frac{1}{m\lambda_{\kstarm}} )$. Combining the last
lower bound, the reduction scheme, and the definition of $\kstarm$
implies the result of the theorem.

To conclude the proof, it remains to check \eqref{pr:theo:lower:m:e1},
$\xdf_\theta\in\Fgr$, and $\edf_\theta\in{\mathcal E}_\lambda
^\edfr
$ for both
$\theta$. To show $\xdf_\theta\in\Fgr$, we first observe that
$\xdf_\theta$ integrates to 1. Moreover, $\xdf_\theta$ is
nonnegative, because $ \vert(1-\theta\alpha_m)
\frac{1\wedge\sqrt{\xdfr-1}}{d^{1/4}} \xdfw_{\kstarm}^{-1/2}\vert
\leq
\xdfw_{\kstarm}^{-1/2}\leq1$ and $\Vert\xdf_\theta\Vert^2_\xdfw=1+
\xdfw_{\kstarm}\vert[\xdf_\theta]_{\kstarm}\vert^2\leq1+
\xdfw_{\kstar}\vert(1-\theta\alpha_m)
\frac{1\wedge\sqrt{\xdfr-1}}{d^{1/4}} \xdfw_{\kstarm}^{-1/2}\vert^2
\leq
\xdfr$. Consider $\edf_\theta$, which obviously integrates to 1.
Furthermore, as $\edf\geq1/2$, the function $\edf_\theta=\edf+
\theta
\alpha_m [\edf]_{\kstarm} e_{\kstarm}$ is nonnegative, because
$\vert\theta
\alpha_m [\edf]_{\kstarm} e_{\kstarm}\vert\leq
\alpha_m\lambda_{\kstarm}^{1/2}\edfr^{1/2} \leq\zeta m^{-1/2}
\sqrt
d\leq1/2$ by using the definition of $\alpha_m$ and $\zeta$. To check
that $\edf_\theta\in{\mathcal E}_\lambda^\edfr$, it remains to
show that
$1/\edfr\leq[\edf_\theta]_j^2/\lambda_j\leq\edfr$ for all
$\vert j\vert>0$. Because $\edf\in{\mathcal E}_\lambda^{\sqrt\edfr
}$, it
follows from the
definition of $\edf_\theta$ that these inequalities are satisfied for
all $j\ne\kstarm$, and, moreover, that $ 1/\edfr\leq\frac{
\vert[\edf]_{\kstarm}\vert^2}{\sqrt{d}\lambda_{\kstarm}}\leq
\frac
{(1+\theta
\alpha_m)^2 \vert[\edf]_{\kstarm}\vert^2}{\lambda_{\kstarm}}\leq
\frac
{\sqrt d
\vert[\edf]_{\kstarm}\vert^2}{\lambda_{\kstarm}} \leq\edfr$.
Finally, consider
\eqref{pr:theo:lower:m:e1}. As in the proof of Theorem~\ref
{theo:lower:n}, by using the Hellinger affinity
$\rho(p_1,p_{-1})$, we obtain, for any estimator $\txdf$ of $\xdf$, that
\[
\bigl\{\Ex_{1}\Vert\txdf-\xdf_{1}\Vert_\hw^2+
\Ex_{-1}\Vert\txdf-\xdf_{1}\Vert_\hw^2
\bigr\} \geq\tfrac{1}{2} \Vert\xdf_{1}-\xdf_{-1}
\Vert^2_\hw\rho(p_1,p_{-1}).
\]
Next, we bound from below the Hellinger affinity $\rho(p_1,p_{-1})\geq
1/4$ for all $m\geq2$, which proves \eqref{pr:theo:lower:m:e1}. From
the independence and the
fact that $\ydf_1=\ydf_{-1}$, it is readily seen that Hellinger
affinity satisfies $\rho(p_{1},p_{-1})= \rho(\ydf_1,\ydf_{-1})^n
\rho(\edf_1,\edf_{-1})^m=\rho(\edf_1,\edf_{-1})^m= (1-\frac
{1}{2}H^2(\edf_{1},\edf_{-1}) )^m$. Thus, we conclude
$\rho(p_{1},p_{-1})\geq(1-1/m)^m \geq1/4$, for all $m\geq2$, because
\begin{eqnarray*}
H^2(\edf_{1},\edf_{-1}) &\leq&\int
\frac{ \vert\edf_{1}-\edf_{-1} \vert^2}{\edf_{1}+\edf_{-1}} =
\int\frac{
\vert\edf_{1}-\edf_{-1} \vert^2}{\edf}\leq2 \int\vert\edf_{1}-
\edf_{-1}\vert^2
\\
&\leq&2 \int4 \alpha_m^2 \bigl\vert[\edf]_{\kstarm}
\bigr\vert^2 \ef_{\kstarm}^2 \leq8 \edfr
\alpha_m^2 \lambda_{\kstarm} = 8\edfr
\zeta^2 m^{-1}\leq2 m^{-1},
\end{eqnarray*}
where we used that $\phi\geq1/2$ and the definition of $\alpha_m$ and
$\zeta$. This completes the proof.
\end{pf*}

\subsubsection*{Upper bound}

\begin{pf*}{Proof of Theorem~\ref{theo:upper}} We begin our proof with the
observation that $\var([\hydf]_j)\leq1/n$ and $\var([\hedf]_j)\leq
1/m$ for all $j\in\Z$. Moreover, by applying Theorem~2.10 of
Petrov \cite{Petrov1995}, there exists a constant $C>0$ such that
$\Ex\vert[\hedf]_j-[\edf]_j\vert^4\leq C/m^2$ for all $j\in\Z$ and
$m\in\N$. We use these results below without further reference.
Now define $\txdf:=1+\sum_{0<\vert j\vert\leq\kstar}
[\xdf]_j\1\{\vert[\hedf]_j\vert^2\geq1/m\} \ef_j$ and decompose
the risk into
two terms,
%
%
\begin{equation}
\label{pr:theo:upper:e1} \Ex\Vert\hxdf_{\kstar}-\xdf\Vert_\hw^2
\leq2\Ex\Vert\hxdf_{\kstar}-\txdf\Vert_\hw^2+2\Ex
\Vert\txdf-\xdf\Vert_\hw^2=: A + B,
\end{equation}
which we bound separately. First, consider $A$, which we decompose
further,
\begin{eqnarray*}
\Ex\Vert\hxdf_{\kstar}- \txdf\Vert_\hw^2 &\leq&2
\sum_{0<\vert j\vert\leq\kstar} \hw_j\E\biggl[
\frac{\vert[\hydf]_j-[\ydf
]_j\vert^2}{\vert[\hedf]_j\vert^2} \mathbh{1} { \bigl\{\bigl\vert
[\hedf]_j
\bigr\vert^2 \geq1/m \bigr\}} \biggr]
\\
&&{}+2\sum_{0<\vert j\vert\leq\kstar}\hw_j \bigl\vert[
\xdf]_j\bigr\vert^2 \E\biggl[\frac{\vert[\hedf]_j-[\edf]_j\vert
^2}{\vert[\hedf]_j\vert^2} \1 \bigl\{
\bigl\vert[ \hedf]_j\bigr\vert^2\geq1/m \bigr\} \biggr]=:
A_1 + A_2.
\end{eqnarray*}
Using the elementary inequality
$ \vert[\edf]_j/[\hedf]_j\vert^2 \leq2 \vert[\edf]_j/[\hedf
]_j-1\vert^2 +2 $,
the independence of $\hedf$ and $\hydf$,
and $\edf\in{\mathcal E}_\lambda^\edfr$, together with the
definition of
$\psi_n$ given in \eqref{def:psin}, we obtain
\[
A_1\leq4 \sum_{0<\vert j\vert\leq\kstar} \hw_j
\biggl\{ \frac{m\var([\hydf]_j)\var([\hedf]_j)}{\vert[\edf
]_j\vert^2} + \frac{\var([\hydf]_j)}{\vert[\edf]_j\vert^2}
\biggr\}\leq8d\sum
_{0<\vert j\vert\leq\kstar} \frac{\hw_j}{n\lambda_j}\leq8d \psi_n.
\]
Moreover, we have $\Ex\frac{\vert[\hedf]_j-[\edf]_j\vert^2}{\vert
[\hedf]_j\vert^2}
\1\{\vert[\hedf]_j\vert^2\geq1/m\}\leq\frac{2m\Ex\vert[\hedf
]_j-[\edf
]_j\vert^4}{\vert[\edf]_j\vert^2} + \frac{2\var([\hedf
]_j)}{\vert[\edf]_j\vert^2}\leq
\frac{ 2(C+1)}{m\vert[\edf]_j\vert^2}\leq\frac{ 2(C+1)\edfr
}{m\lambda_j}$
and $\Ex\frac{\vert[\hedf]_j-[\edf]_j\vert^2}{\vert[\hedf
]_j\vert^2}
\1\{\vert[\hedf]_j\vert^2\geq1/m\}\leq1$, where we have again used the
elementary inequality and $\edf\in{\mathcal E}_\lambda^\edfr$.
Combining both
bounds together with $\xdf\in\Fgr$ and the definition of $\kappa_m$
given in \eqref{def:kappam}, we obtain
\[
A_2\leq4 (C+1) \edfr\sum_{0<\vert j\vert\leq\kstar}
\hw_j\bigl\vert[\xdf]_j\bigr\vert^2 \min\biggl(1,
\frac{1}{m\lambda_j} \biggr)\leq4(C+1) \edfr\xdfr\kappa_{m}.
\]
Now consider $B$, which we decompose further into
\begin{eqnarray*}
\Ex\Vert\txdf- \xdf\Vert_\hw^2&=&\sum
_{0<\vert j\vert} \hw_j\bigl\vert[\xdf]_j
\bigr\vert^2 \bigl(1- \1\bigl\{0<\vert j\vert\leq\kstar\bigr\}\1 \bigl\{\bigl\vert[
\hedf]_j\bigr\vert^2 \geq1/m \bigr\} \bigr)^2
\\
&=& \sum_{\vert j\vert> \kstar} \hw_j\bigl\vert[
\xdf]_j\bigr\vert^2 +\sum_{0<\vert j\vert\leq\kstar}
\hw_j \bigl\vert[\xdf]_j\bigr\vert^2 \P\bigl(\bigl\vert[
\hedf]_j\bigr\vert^2<1/m \bigr)=: B_1 +
B_2,
\end{eqnarray*}
where $B_1\leq\Vert\xdf\Vert^2_\xdfw\hw_{\kstar}\xdfw_{\kstar
}^{-1}\leq\xdfr\psi_n$, because $\xdf\in\Fgr$. Moreover, $B_2\leq
4 \edfr\xdfr\kappa_{m}$, using that
%
%
\begin{equation}
\label{pr:theo:upper:e2} \P\bigl(\bigl\vert[\hedf]_j
\bigr\vert^2<1/m \bigr)\leq4\edfr\min\biggl(1,\frac{1}{m \lambda_j}
\biggr),
\end{equation}
which we show below. The result of the theorem now follows by combining
the decomposition \eqref{pr:theo:upper:e1} and the estimates of $A_1,
A_2, B_1$, and $B_2$.

To conclude, we prove \eqref{pr:theo:upper:e2}. If $\vert[\edf
]_j\vert^2\geq
4/m$, then, using Tchebychev's inequality, we deduce that 
\begin{eqnarray*}
\P\bigl(\bigl\vert\hfedf_j\bigr\vert^2<1/m \bigr)&\leq&\P\bigl(\bigl\vert
\hfedf_j/[\edf]_j\bigr\vert<1/2\bigr) \leq\P\bigl(\bigl\vert
\hfedf_j- [\edf]_j\bigr\vert>\bigl\vert[\edf]_j\bigr\vert
/2\bigr)
\\
&\leq&4 \frac{\var
(\hfedf_j)}{\vert[\edf]_j\vert^2} \leq4\edfr/(m \lambda_j).
\end{eqnarray*}\vspace*{-1pt}
On the other hand, in the case where $\vert[\edf]_j\vert^2<4/m$, the estimate
$\P(\vert\hfedf_j\vert^2<1/m)\leq4\edfr/(m \lambda_j)$ also
holds, because
$1\leq
4/(m\vert[\edf]_j\vert^2)\leq4\edfr/(m \lambda_j) $. Combining
the last
estimates and $\P(\vert\hfedf_j\vert^2<1/m)\leq1$, we obtain
\eqref
{pr:theo:upper:e2}, which completes the proof.
\end{pf*}

\subsubsection*{Illustration: Estimation of derivatives}

\begin{pf*}{Proof of Proposition~\ref{coro:ex:lower}} Because for each
$0\leq s
\leq p$, we have $\Ex\Vert\txdf^{(s)}-\xdf^{(s)}\Vert^2 \sim
\Ex
\Vert\txdf-\xdf\Vert_{\hw}^2$, we intend to apply the general result
given in Corollary~\ref{cor:lower}. In both cases, the additional
conditions formulated in Theorem~\ref{theo:lower:n} and~\ref
{theo:lower:m} are readily verified. Thus, it is sufficient to
evaluate the lower bounds $\psi_n$ and $\kappa_m$ given in
\eqref{def:psin} and \eqref{def:kappam}, respectively. Note that the
optimal dimension parameter,
$\kstar:=\argmin_{j\in\N}[\max(\frac{\hw_j}{\xdfw_j}, \sum
_{0<\vert l\vert\leq
j}\frac{\hw_l}{n\lambda_l})]$ satisfies $n\hw_{\kstar}/\xdfw
_{\kstar}
\sim\sum_{0<\vert l\vert\leq{\kstar}}\hw_l/\lambda_l$, because
both sequences
$(\xdfw_j/\hw_j)$ and $(\sum_{0<\vert l\vert\leq j}\frac{\hw
_l}{n\lambda_l})$ are
non-increasing.

{[o-o]} The well-known approximation $\sum_{j=1}^{m} j^{r}\sim
m^{r+1}$ for $r>0$ implies that
$(\xdfw_{\kstar}/\hw_{\kstar})\*\sum_{0<\vert l\vert\leq{\kstar
}}\hw_l/\lambda_l
\sim(\kstar)^{2a+2p+1}$. It follows that $\kstar\sim n^{1/(2p+2a+1)}$,
and the first lower bound is $\psi_n \sim
n^{-(2p-2s)/(2p+2a+1)}$. Moreover, we have $\kappa_m\sim
m^{-([p-s]\wedge a)/a} $, because the minimum in $\kappa_m=
\sup_{j\in\Z}\{\vert j\vert^{-2(p-s)} \min(1,\vert j\vert
^{2a}/m)\}$ is equal to 1
for $\vert j\vert\geq m^{1/2a}$ and $\vert j\vert^{-2(p-s)}$ is
non-increasing.

{[s-o]} Approximating the sum in the same way as above, we obtain
$(\xdfw_{\kstar}/\hw_{\kstar})\sum_{0<\vert l\vert\leq{\kstar
}}\hw_l/\allowbreak\lambda_l
\sim(\kstar)^{2a+1}\exp({\kstar}^{2p})$, and thus $\kstar\sim
(\log
n)^{1/(2p)}$. The resulting\vspace*{1pt} rate is $\psi_n \sim n^{-1} \*(\log
n)^{(2a+2s+1)/(2p)}$. Furthermore, we have $\kappa_m\sim m^{-1}$,
because the supremum is taken over
$j^{2s}\exp(-j^{2p}) \min(1,j^{2a}/m)$, which takes its maximum
at the border because of the dominating exponential term.

{[o-s]} Applying Laplace's method (see chap. 3.7 in
Olver \cite{Olver1974}), we have
$(\xdfw_{\kstar}/\hw_{\kstar})\*\sum_{0<\vert l\vert\leq
{\kstar}}\hw_l/\lambda_l\sim(\kstar)^{2p + ((2a-1)\vee
0)}\exp(\vert\kstar\vert^{2a})$, which implies that $\kstar\sim
(\log
n)^{1/(2a)}$ and that the first lower bound can be rewritten as
$\psi_n\sim(\log n)^{-(p-s)/a}$. Furthermore, we have $\kappa_m\sim
(\log m)^{-(p-s)/a}$, because the minimum in $\kappa_m=
\sup_{j\in\Z}\{\vert j\vert^{-2(p-s)} \min(1,\break\exp(\vert j\vert
^{2a})/ m)\}$ is equal to
1 for $\vert j\vert\geq(\log m)^{(1/2a)}$ and $\vert j\vert
^{-2(p-s)}$ is
non-increasing. Consequently, the lower bounds in Proposition~\ref
{coro:ex:lower} follow by applying Corollary~\ref{cor:lower}.
\end{pf*}

\begin{pf*}{Proof of Proposition~\ref{coro:ex:upper}} The result is an
immediate
consequence of Theorem~\ref{theo:upper} and Proposition~\ref{coro:ex:lower}.
\end{pf*}

\subsection{\texorpdfstring{Proofs of Section~\protect\ref{sec:adaptive-known}}
{Proofs of Section 3}}\label{app:proofs:known}

\subsubsection*{Partially adaptive estimation}
\label{sec:part-adapt-estim}

We begin by defining and recalling notations to be used in the
proof. Given $u\in L^2[0,1]$, we denote by $[u]$ the infinite vector of
Fourier coefficients $[u]_j:=\langle u,e_j\rangle$. In particular, we use
the notations
\begin{eqnarray*}
&&\hspace*{5.7pt}\hxdf_k= \sum_{j=- k}^k
\frac{\widehat{[\ydf]}_j}{\widehat{[\edf]}_j}\1 \bigl\{\bigl\vert
\widehat{[\edf ]}_j\bigr\vert^2
\geq1/m \bigr\} e_j,\qquad\txdf_k:= \sum
_{j=- k}^k\frac{\widehat{[\ydf]}_j }{[\edf]_j}e_j,\qquad
\xdf_k:=\sum_{j=- k}^k
\frac{[\ydf]_j}{
[\edf]_j} e_j,
\\
&&\hPhi_{u}:= \sum_{j\in\Z}
\frac{[u]_j}{\widehat{[\edf]}_j}\1 \bigl\{\bigl\vert\widehat{[\edf
]}_j\bigr\vert^2
\geq1/m \bigr\}e_j,\qquad\tPhi_{u}:= \sum
_{j\in\Z}\frac{[u]_j}{[\edf]_j}e_j.
\end{eqnarray*}
Furthermore, let $\hydf$ be the function with Fourier coefficients
$[\hydf]_j:=\widehat{[\ydf]}_j$. Given $0\leq k\leq k'$, we then have,
for all $t\in\cS_k:=\operatorname{span}\{e_{-k},\ldots,e_k\}$,
\begin{eqnarray*}
\langle t, \xdf_{k'}\rangle_\hw&=& \langle t,
\tPhi_{\ydf}\rangle_\hw= \sum_{j=- k}^k
\frac{\hw_j[t]_j[\ydf]_j}{[\edf]_j}= \sum_{j=- k}^k
\hw_j[t]_j[\xdf]_j=\langle t, \xdf
\rangle_\hw,
\\
\langle t, \txdf_{k'}\rangle_\hw&=&\langle t,
\tPhi_{\hydf}\rangle_\hw= \frac{1}{n} \sum
_{i=1}^n \sum_{j=-k}^ke_j(-Y_i)
\frac{\hw_j[t]_j}{[\edf]_j}= \langle t, \txdf_{k}\rangle_\hw,
\\
\langle t, \hxdf_{k'}\rangle_\hw&=& \langle t,
\hPhi_{\hydf}\rangle_\hw=\frac{1}{n} \sum
_{i=1}^n \sum_{j=-k}^ke_j(-Y_i)
\frac{\hw_j[t]_j}{\widehat{[\edf]}_j}\1 \bigl\{\bigl\vert\widehat
{[\edf ]}_j\bigr\vert^2
\geq1/m \bigr\}=\langle t, \hxdf_{k}\rangle_\hw.
\end{eqnarray*}
Consider the function $\nu=\hydf-\ydf$ with Fourier coefficients
$[\nu]_j =\widehat{[\ydf]}_j - [\ydf]_j = \widehat{[\ydf]}_j -
\Ex\widehat{[\ydf]}_j$. We then have, for every $t\in\cS_k$,
%
%
\begin{eqnarray}
\label{pr:theo:adap1:def2} \langle t,\hPhi_{\hg} - f
\rangle_\hw&=& \langle t,\hPhi_{\hg} - \tPhi_g
\rangle_\hw= \langle t,\tPhi_{\hg} - \tPhi_g
\rangle_\hw+ \langle t,\hPhi_{\hg} - \tPhi_{\hg}
\rangle_\hw
\nonumber
\\[-8pt]
\\[-8pt]
&=& \langle t,\tPhi_{\nu}\rangle_\hw+ \langle t,
\hPhi_{\hg} - \tPhi_{\hg}\rangle_\hw=\langle t,
\tPhi_{\nu}\rangle_\hw+\langle t,\hPhi_{\nu} -
\tPhi_{\nu}\rangle_\hw+\langle t,\hPhi_{\ydf} -
\tPhi_{\ydf}\rangle_\hw.
\nonumber
\hspace*{27pt}
\end{eqnarray}
At the end of this section, we prove some technical lemmas that
are used in the following proof.

\begin{pf*}{Proof of Theorem~\ref{theo:phi}}
We consider the contrast
\begin{eqnarray*}
\ct(t) := \Vert t\Vert^2_\hw- 2 {\langle t,
\hPhi_{\hydf}\rangle}_\hw\qquad\forall t\in
L^2[0,1].
\end{eqnarray*}
It obviously follows that, for all $t\in\cS_k$,
$\ct(t) = \Vert t-\hxdf_k\Vert^2_\hw- \Vert\hxdf_k\Vert^2_\hw$,
and thus,
%
%
\begin{equation}
\label{sec:gen:con:2} \arg\min_{t\in\cS_k} \ct(t) = \hxdf_{k}
\qquad\forall k\geq0.
\end{equation}
Moreover, the adaptive choice of the dimension parameter from \eqref
{gen:ada:1} can be rewritten as
%
%
\begin{equation}
\label{sec:gen:ada:1} \tk= \mathop{\argmin}_{0\leq k\leq(\Cy
_n\wedge\Ce
_m)} \biggl[\ct(
\hxdf_{k}) + 60\frac{\delta_k}{n} \biggr].
\end{equation}
Let $\pen(k):= 60 \delta_k/n$; then, for all $1\leq k \leq \Mnm$,
we have
\begin{eqnarray*}
\ct(\hxdf_{\tk}) + \pen(\tk) \leq\ct(\hxdf_{k}) + \pen(k)
\leq\ct(\xdf_k) + \pen(k),
\end{eqnarray*}
using first \eqref{sec:gen:ada:1} and then \eqref{sec:gen:con:2}.
This inequality
implies that
\begin{eqnarray*}
\Vert\hxdf_{\tk}\Vert^2_\hw-\Vert
\xdf_k\Vert^2_\hw&\leq2\langle
\hxdf_{\tk}-\xdf_k, \hPhi_{\hydf}
\rangle_\hw+ \pen(k)-\pen(\tk),
\end{eqnarray*}
and thus, using \eqref{pr:theo:adap1:def2}, we have, for all $1\leq k
\leq \Mnm$,
%
%
\begin{eqnarray}
\label{eq:5:1} %
\Vert\hxdf_{\tk}-\xdf\Vert^2_\hw
&\leq&\Vert\xdf-\xdf_{k}\Vert^2_\hw+ \pen(k)-
\pen(\tk)
\nonumber
\\[-8pt]
\\[-8pt]
&&{}+ 2\langle\hxdf_{\tk} - \xdf_k, \tPhi_\nu
\rangle_\hw+ 2\langle\hxdf_{\tk} - \xdf_k,
\hPhi_{\nu} - \tPhi_{\nu}\rangle_\hw+2\langle
\hxdf_{\tk} - \xdf_k,\hPhi_{\ydf} -
\tPhi_{\ydf
}\rangle_\hw.
\nonumber
\hspace*{40pt} %
\end{eqnarray}
Consider the unit ball $\cB_k:=\{f\in\cS_k\dvt\Vert f\Vert_\hw
\leq1\}$
and, for arbitrary $\tau>0$ and $t\in\cS_k$, the elementary
inequality
\[
2\bigl\vert\langle t,h\rangle_\hw\bigr\vert\leq2\Vert t\Vert_\hw
\sup_{t\in\cB_k}\bigl\vert\langle t,h\rangle_\hw\bigr\vert\leq\tau\Vert
t \Vert_\hw^2+\frac{1}{\tau} \sup_{t\in\cB_k}
\bigl\vert\langle t,h\rangle_\hw\bigr\vert^2= \tau\Vert t
\Vert_\hw^2+ \frac{1}{\tau}\sum
_{j=-k}^k \hw_j\bigl\vert
[h]_j\bigr\vert^2.
\]
Combining the last estimate with (\ref{eq:5:1}) and $\hxdf_{\tk
}-\xdf_k \in\cS_{\tk\vee k} \subset\cS_{\Cy_n\wedge\Ce_m}$, we obtain
\begin{eqnarray*}
\Vert\hxdf_{\tk} -\xdf\Vert^2_\hw&\leq&\Vert
\xdf-\xdf_k\Vert^2_\hw+ 3\tau\Vert
\hxdf_{\tk} -\xdf_k\Vert^2_\hw+
\pen(k) - \pen(\tk)
\\
&&{}+ \frac{1}{\tau}\supt{k\vee\tk}\bigl\vert\langle t,\tPhi_\nu
\rangle_\hw\bigr\vert^2 + \frac{1}{\tau}\supt{\Mnm}\bigl\vert
\langle t, \hPhi_{\nu} - \tPhi_{\nu
}\rangle_\hw
\bigr\vert^2
\\
&&{}+ \frac{1}{\tau}\supt{ \Mnm}\bigl\vert\langle t,\hPhi_{\ydf} -
\tPhi_{\ydf}\rangle_\hw\bigr\vert^2.
\end{eqnarray*}
Note that $\Vert\hxdf_{\tk}-\xdf_k\Vert^2_\hw
\leq2\Vert\hxdf_{\tk}-\xdf\Vert^2_\hw+2\Vert\xdf_k-\xdf\Vert
^2_\hw
$ and
that $\Vert\xdf-\xdf_k\Vert^2_\hw\leq\xdfr\hw_k/\xdfw_k$ for all
$\xdf\in\Fgr$, because $\hw/\xdfw$ is non-increasing. Setting
$\tau:=
{1}/{8}$, we obtain
%
%
\begin{eqnarray}
\label{tau:acht} \frac{1}{4}\Vert\hxdf_{\tk} - \xdf
\Vert^2_\hw&\leq&\frac{7}{4} (\xdfr
\hw_k / \xdfw_k) + \pen(k) - \pen(\tk)
\nonumber
\\
&&{}+ 8\supt{k\vee\tk}\bigl\vert\langle t,\tPhi_\nu
\rangle_\hw\bigr\vert^2 + 8 \supt{\Mnm}\bigl\vert\langle t,
\hPhi_{\nu} - \tPhi_{\nu}\rangle_\hw
\bigr\vert^2
\\
&&{}+ 8\supt{\Mnm}\bigl\vert\langle t,\hPhi_{\ydf} - \tPhi_{\ydf}
\rangle_\hw\bigr\vert^2.
\nonumber
\end{eqnarray}

Defining the event
%
%
\begin{equation}
\label{eq:20} \Omega_q := \biggl\{ \forall0\leq\vert j\vert\leq
M_m^u \dvt\biggl\vert\frac{1}{\widehat{[\edf]}_j}- \frac{1}{[\edf]_j}
\biggr\vert\leq\frac{1}{2\vert[\edf]_j\vert} \wedge\bigl\vert\widehat
{[\edf]}_j
\bigr\vert^2 \geq1/m \biggr\},
\end{equation}
consider the following decomposition of the risk:
%
%
\begin{equation}
\label{eq:risiko:q} \E\Vert\hxdf_{\tk} -\xdf\Vert^2_\hw
= \E\Vert\hxdf_{\tk} -\xdf\Vert^2_\hw\mathbh{1}
{\{\Omega_q\}} + \E\Vert\hxdf_{\tk} -\xdf
\Vert^2_\hw\mathbh{1} { \bigl\{\Omega_q^c
\bigr\}}.
\end{equation}
We bound these two terms separately.
Consider the first term. By Lemma~\ref{lem:schachtel} below and
$\mathbh{1}{\{\vert\widehat{[\edf]}_j\vert^2\geq1/m\}}\mathbh
{1}{\{
\Omega_q\}}
= \mathbh{1}{\{\Omega_q\}}$,
it follows that for all $1\leq\vert j\vert\leq\Mnm$,
\begin{eqnarray*}
\biggl(\frac{[\edf]_j}{\widehat{[\edf]}_j}\mathbh{1} { \bigl\{
\bigl\vert
\widehat{[\edf]}_j
\bigr\vert^2 \geq1/m \bigr\}} - 1 \biggr)^2 \mathbh{1} {\{
\Omega_q\}} = \bigl\vert[\edf]_j\bigr\vert^2
\mathbh{1} {\{\Omega_q \}} \biggl\vert\frac{1}{\widehat{[\edf]}_j}-
\frac{1}{[\edf]_j}\biggr \vert^2 \leq\frac{1}{4}.
\end{eqnarray*}
Thus,
$ \supt{k} \vert\langle t,\hPhi_{\nu} - \tPhi_{\nu}\rangle_\hw
\vert^2
\mathbh{1}{\{\Omega_q\}}
\leq\frac{1}{4} \supt{k} \vert\langle t,\tPhi_\nu\rangle_\hw
\vert^2$
for all $0\leq k\leq\Mnm$, and \eqref{tau:acht} implies that
%
%
\begin{eqnarray}\label{eq:12}
&&\frac{1}{4} \Vert\hxdf_{\tk} - \xdf
\Vert^2_\hw\mathbh{1} {\{\Omega_q\}}\nonumber\\
&&\quad{}\leq\frac{7}{4} (\xdfr\hw_k / \xdfw_k) + 10 \Bigl(
\supt{k\vee\tk} \bigl\vert\langle t,\tPhi_\nu\rangle_\hw
\bigr\vert^2 - (6 \delta_{k\vee\tk})/n \Bigr)_+\\
&&\qquad{}+ (60 \delta_{k\vee\tk} )/n + \pen(k) - \pen(\tk) + 8\supt
{\Mnm}\bigl\vert\langle t,\hPhi_{\ydf} - \tPhi_{\ydf}\rangle_\hw\bigr\vert^2.
\nonumber
\end{eqnarray}
Moreover, we have that $ 60 \delta_{k\vee
\tk}/n=\pen(k\vee\tk)\leq\pen(k)+\pen(\tk)$. Further note that
%
%
\begin{equation}
\label{eq:deltaDeltazeta} \Delta_k \leq\edfr\Delta_k^\lambda,
\qquad\delta_k\leq\edfr\zeta_\edfr\delta_k^\lambda,
\quad\mbox{and}\quad\delta_k/ \Delta_k \geq2 k
\zeta_\edfr^{-1} \frac{\log(\Delta_k^\lambda\vee(k+2))}{\log(k+2)},
\end{equation}
with $\zeta_\edfr= \log(3\edfr)/\log\edfr$.
From Lemma~\ref{lem:schachtel}, it follows that
\begin{eqnarray*}
&&\sup_{\xdf\in\Fgr}\sup_{\edf\in\Eld}\E\Vert\hxdf_{\tk}
-\xdf\Vert_\hw^2\mathbh{1} {\{\Omega_q\}}
\\
&&\quad\leq480 (\xdfr+ \edfr\zeta_\edfr) \min_{0\leq k\leq\Cy
_n^\lambda
\wedge\Ce_m^\lambda} \bigl[ \max
\bigl( {\hw_k}/{\xdfw_k} , {\delta_k^\lambda}/{n}
\bigr) \bigr]
\\
&&\qquad{}+ 40 \sup_{\xdf\in\Fgr}\sup_{\edf\in\Eld} \sum
_{0\leq k'\leq(\Cy_n^u\wedge\Ce_m^u)}\E\Bigl( \supt{k'} \bigl\vert
\langle t,
\tPhi_\nu\rangle_\hw\bigr\vert^2 - (6
\delta_{k'})/n \Bigr)_+
\\
&&\qquad{}+ 32 \sup_{\xdf\in\Fgr}\sup_{\edf\in\Eld}\E\Bigl[\supt{
\bigl(
\Cy_n^u\wedge\Ce^u_m \bigr)}\bigl\vert
\langle t,\hPhi_{\ydf} - \tPhi_{\ydf}\rangle_\hw
\bigr\vert^2 \Bigr].
\end{eqnarray*}
To bound the second term, we apply Lemma~\ref{lem:talalem}
with $\delta_k^* = \delta_k$ and $\Delta_k^* =
\Delta_k$. By virtue of \eqref{eq:deltaDeltazeta}, we have, for all
$k\geq0$,
\begin{eqnarray*}
&&\E\biggl( \supt{k} \bigl\vert\langle t,\tPhi_\nu\rangle_\hw
\bigr\vert^2 - 6 \frac{\delta_k}{n} \biggr)_+
\\
&&\quad\leq C \biggl\{ \frac{1}{n^2} \exp\bigl( -K_2 \sqrt{n} \bigr) \edfr
\zeta_\edfr\delta_k^\lambda
\\
&&\hphantom{\quad\leq C \biggl\{} {}+\frac{\Vert\edf\Vert^2
\Vert\xdf\Vert^2}{n} \edfr
\Delta_k^\lambda\exp\biggl(- \frac{k}{3 \Vert\edf\Vert^2 \Vert
\xdf\Vert^2 \zeta_\edfr}
\frac
{\log(\Delta_k^\lambda\vee(k+2))}{\log(k+2)} \biggr) \biggr\}.
\end{eqnarray*}
Owing to Lemmas~\ref{lem:schachtel} and~\ref{lem:NuMu}(i) and the
properties of the function
$\Sigma$ from Definition~\ref{def:known}, we have
\begin{eqnarray*}
\sum_{k=0}^{\Cy_n^u} \E\biggl( \supt{k} \bigl\vert
\langle t,\tPhi_\nu\rangle_\hw\bigr\vert^2 - 6
\frac{\delta_k}{n} \biggr)_+ \leq\frac{C}{n} \edfr\Sigma\bigl
(\Vert\edf
\Vert^2 \Vert\xdf\Vert^2 \zeta_\edfr\bigr).
\end{eqnarray*}
It can be readily verified that $\Vert\edf\Vert^2\leq\edfr\Lambda
$ for
all $\edf\in\Eld$ and $\Vert\xdf\Vert^2\leq\xdfr$ for all
$\xdf\in\Fgr$. The remaining term can be controlled by virtue of
Lemma~\ref{lem:Q}, which shows that
%
%
\begin{eqnarray}
\label{deco:known:a} \sup_{\xdf\in\Fgr}\sup_{\edf\in\Eld} \E
\Vert
\hxdf_{\tk} -\xdf\Vert^2_\hw\mathbh{1} {\{
\Omega_q\}} &\leq& C \Bigl\{ (\xdfr+ \edfr\zeta_\edfr)
\min_{0\leq k\leq(\Cy
_n^\lambda\wedge
\Ce_m^\lambda)} \bigl[ \max\bigl({\hw_k}/{\xdfw_{k}},
{\delta_k^\lambda}/{n} \bigr) \bigr]
\nonumber
\\[-8pt]
\\[-8pt]
&&\hphantom{C \bigl\{} {}+ \xdfr\edfr\kappa_m + \edfr\Sigma
(\xdfr
\edfr\Edfw\zeta_\edfr) n^{-1} \Bigr\}.
\nonumber
\end{eqnarray}

Consider the second term from \eqref{eq:risiko:q}.
Let $\breve\xdf_{k}:=1+\sum_{0<\vert j\vert\leq
k} [\xdf]_j \1\{\vert\hfedf_j\vert^2\geq1/m\}\ef_j$. It is easy to
see that $\Vert\hxdf_{k}-\breve\xdf_{k}\Vert^2 \leq
\Vert\hxdf_{k'}-\breve\xdf_{k'}\Vert^2$ for all $ k\leq k'$ and
$\Vert\breve\xdf_{k}-\xdf\Vert^2\leq\Vert\xdf\Vert^2$ for all
$k\geq
0$. Thus, using that $0\leq\tk\leq({\Cy_n^\circ}\wedge m)$, we
can write
\begin{eqnarray*}
\Ex\Vert\hxdf_{\tk}-\xdf\Vert_\omega^2\mathbh{1}
{ \bigl\{\Omega_{q}^c \bigr\}} &\leq&2 \bigl\{\Ex\Vert
\hxdf_{\tk}-\breve\xdf_{\tk
}\Vert_\omega^2
\mathbh{1} { \bigl\{\Omega_{q}^c \bigr\}} + \Ex\Vert
\breve\xdf_{\tk}- \xdf\Vert_\omega^2\mathbh{1} {
\bigl\{ \Omega_{q}^c \bigr\}} \bigr\}
\\
&\leq&2 \bigl\{ \Ex\Vert\hxdf_{({\Cy_n^\circ}\wedge m)}-\breve
\xdf_{({\Cy_n^\circ}\wedge m)}
\Vert_\omega^2 \mathbh{1} { \bigl\{\Omega_{q}^c
\bigr\}} + \Vert\xdf\Vert_\omega^2 \P\bigl[
\Omega_q^c \bigr] \bigr\}.
\end{eqnarray*}
Moreover, applying Theorem~2.10 of Petrov \cite{Petrov1995},
\begin{eqnarray*}
&&\Ex\Vert\hxdf_{({\Cy_n^\circ}\wedge m)}-\breve\xdf_{({\Cy
_n^\circ}\wedge m)}\Vert_\omega^2
\mathbh{1} { \bigl\{\Omega_{q}^c \bigr\}}
\\
&&\quad\leq2m\sum_{0<\vert j\vert\leq({\Cy_n^\circ}\wedge
m)}\omega_j \bigl
\{ \Ex\bigl( \hfydf_j -[\edf]_j [\xdf]_j\bigr)^2
\mathbh{1} { \bigl\{\Omega_{q}^c \bigr\}} + \Ex\bigl(
[\edf]_j [\xdf]_j-\hfedf_j
[\xdf]_j\bigr)^2\mathbh{1} { \bigl\{\Omega_{q}^c
\bigr\}} \bigr\}
\\
&&\quad\leq2 m \biggl\{ \sum_{0<\vert j\vert\leq({\Cy_n^\circ}\wedge
m)} \omega_j
\bigl[\Ex\bigl(\hfydf_j - [\ydf]_j \bigr)^{4}
\bigr]^{1/2} \P\bigl[\Omega_{q}^c
\bigr]^{1/2}
\\
&&\hphantom{\quad\leq2 m \biggl\{} {}+ \sum_{0<\vert j\vert\leq
({\Cy_n^\circ}\wedge m)}
\omega_j \bigr\vert[\xdf]_{j}\bigl\vert^2 \bigl[\Ex\bigl(
\hfedf_j- [\edf]_j\bigr)^4 \bigr]^{1/2}
\P\bigl[ \Omega_{q}^c \bigr]^{1/2} \biggr\}
\\
&&\quad\leq2 m \Bigl\{ 2 m \Bigl({\max_{1\leq j \leq{\Cy_n^\circ
}}\omega_j}
\Bigr) \bigl(Cn^{-1} \bigr) + \bigl(Cm^{-1} \bigr) \Vert\xdf
\Vert_\omega^2 \Bigr\} \P\bigl[\Omega_{q}^c
\bigr]^{1/2},
\end{eqnarray*}
which implies, using Definition~\ref{def:known}(ii),
%
%
\begin{eqnarray}
\label{eq:5} %
\Ex\Vert\hxdf_{\tk}-\xdf\Vert_\omega^2
\mathbh{1} { \bigl\{\Omega_{q}^c \bigr\}} &\leq&4C
\bigl( m^2 + \Vert\xdf\Vert_\omega^2 \bigr) \P
\bigl[ \Omega_{q}^c \bigr]^{1/2} + 2\Vert\xdf
\Vert_\omega^2 \P\bigl[ \Omega_q^c
\bigr]
\nonumber
\\[-8pt]
\\[-8pt]
& \leq&6Cm^2 \bigl(1+\Vert f\Vert_\hw^2
\bigr)\P\bigl[\Omega_q^c \bigr]^{1/2}.
\nonumber
\end{eqnarray}
By Lemma~\ref{lem:POmegaqc}, it follows that for all $m\in\N$,
%
%
\begin{equation}
\label{eq:4} \sup_{\xdf\in\Fgr}\sup_{\edf\in\Eld} \Ex\Vert
\hxdf_{\whk}-\xdf\Vert_\omega^2\mathbh{1} { \bigl\{
\Omega_{p}^c \bigr\}} \leq C(\edfr) (1+\xdfr)
m^{-1}.
\end{equation}
The result of the theorem
follows from a combination of the last estimate and \eqref{deco:known:a}.
\end{pf*}

\setcounter{lem}{0}
\renewcommand{\thelem}{A\arabic{lem}}
\begin{lem}\label{lem:schachtel}
Under Assumption~\textup{\ref{ass:minreg}}, we have, for all $n,m\in
\N$,
\[
\Cy^\lambda_n \leq\Cy_n \leq
\Cy^u_n \quad\mbox{and}\quad\Ce^\lambda_m
\leq\Ce_m \leq\Ce^u_m.
\]
\end{lem}

\begin{pf}
We first prove that $\Cy_n^\lambda\leq\Cy_n$. If $\Cy
_n^\lambda=
0$ or $\Cy_n = \Cy_n^\circ$, then there is nothing to show. Noting that
\begin{eqnarray*}
\Cy_n^\lambda= 0 &\quad\iff\quad&\max_{1\leq j \leq\Cy_n^\circ}
\frac{\lambda_j}{j\hw_j^+} < \frac{4\edfr\log n}{n} \quad\mbox{and}\\
\Cy_n = 0 &\quad\iff\quad&\max_{1\leq j \leq\Cy_n^\circ} \frac{\lambda_j}{j\hw_j^+} <
\frac{\edfr\log n}{n},
\end{eqnarray*}
we deduce that in the case where $\Cy_n=0$, we also have $\Cy
_n^\lambda
= 0$. This also holds when $\Cy_n^\lambda> 0$ and $\Cy_n^\circ>\Cy_n
> 0$, which
implies
\[
\min_{1\leq j \leq\Cy_n^\lambda} \frac{\lambda_j}{j\hw_j^+} \geq
\frac{4\edfr\log n}{n}\quad\mbox{and}
\quad\frac{\log n}{n} > \frac{\vert[\edf]_{\Cy_n+1}\vert^2}{\Cy
_n \hw_{\Cy_n+1}} \geq\frac{\lambda_{\Cy_n+1}}{\edfr\Cy_n\hw
_{\Cy_n+1}^+}
\]
and thus $\Cy_n+1 > \Cy_n^\lambda$, which proves the claim.

We now prove $\Cy_n\leq\Cy_n^u$. If $\Cy_n = 0$ or $\Cy_n^u = n$,
then this is trivial. On the other hand, if $n>\Cy_n^u\geq0$ and
$\Cy_n^\circ\geq\Cy_n>0$, then it follows from the definitions that
\[
\min_{1\leq j \leq\Cy_n} \frac{\edfr\lambda_j}{j\hw_j^+} \geq
\min_{1\leq
j\leq\Cy_n}
\frac{\vert[\edf]_j\vert^2}{j\hw_j^+} \geq\frac{\log n}{n}\quad
\mbox{and}\quad
\frac{\lambda_{\Cy_n^\circ+1}}{(\Cy_n^\circ+ 1)\hw_{\Cy_n^\circ
+1}^+}< \frac{\log n}{4\edfr n},
\]
which implies that $\Cy_n^\circ+1>\Cy_n$, and hence the claim.
Similar arguments show the corresponding estimates in $m$.
\end{pf}

\begin{lem}\label{lem:NuMu}
Under Assumption~\textup{\ref{ass:minreg}}, we have that
\begin{longlist}[(ii)]
\item
$ \delta_{\Cy^u_n} / n \leq32 \edfr^2 $ for all $n\geq1$,
\item$ m^7 \exp( - \frac{m \lambda_{\Ce^u_m}}{72 \edfr}
) \leq C(\edfr)$ for all $m\geq1$,
\end{longlist}
and for $m\geq\exp(512\log(3\edfr)^2)$ that
\begin{longlist}[(iii)]
\item[(iii)]
$\min_{1\leq j \leq\Ce^u_m} \vert[\edf]_j\vert^2\geq
\frac{2}{m}$.
\end{longlist}
\end{lem}
\begin{pf}
(i) For $\Cy_n^u = 0$, we have $\delta_{\Cy_n^u} = 0$, and there is
nothing to show. If $0<\Cy_n^u\leq n$, then we can show that
$\hw_{\Cy_n^u}^+/\lambda_{\Cy_n^u}\leq4\edfr n / (\Cy_n^u\log(n+2))$,
which we use in the following computation:
\begin{eqnarray*}
\delta_{\Cy_n^u} &=& \Cy_n^u \frac{\hw_{\Cy_n^u}^+}{\lambda_{\Cy_n^u}}
\frac{\log((\hw_{\Cy_n^u}^+/\lambda_{\Cy_n^u}) \vee(\Cy_n^u + 2)
)}{\log(\Cy_n^u + 2)}
\\
&\leq&\frac{4\edfr n}{\log(n+2)} {\log\biggl( \frac{4\edfr n}{\Cy
_n^u\log(n+2)} \vee\bigl(
\Cy_n^u + 2\bigr) \biggr)}\Big/{\log\bigl(
\Cy_n^u + 2\bigr)}
\\
&\leq& n %
\cases{ 4\edfr&\quad$ \bigl(\log(n+2)\geq4\edfr\bigr)$,
\cr
4 \edfr\bigl(4\edfr+ \log(4\edfr) \bigr)/ \bigl(\log(n+2) \bigr
) &\quad
(otherwise),} %
\end{eqnarray*}
which implies that $\delta_{\Cy_n^u}/n \leq4\edfr(4\edfr+
\log(4\edfr))\leq32\edfr^2$
for all $n\geq1$.

(ii) For $0<\Ce_m^u\leq m$, we have
$\lambda_{\Ce_m^u} \geq m^{-1+b_m} (4\edfr)^{-1}$. Thus,
\[
m^7\exp\biggl(-\frac{m\lambda_{\Ce_m^u}}{72\edfr} \biggr) \leq
\exp\biggl(-
\frac{m^{b_m}}{288\edfr^2} + 7\log m \biggr).
\]
This proves the claim, because $\log m \lesssim m^{b_m}$. Note that
$\Ce_m^u = 0$ cannot occur, because we assume that $\lambda_1=1$.

(iii) We have that
\[
\min_{1\leq j\leq\Ce_m^u} \bigl\vert[\edf]_j\bigr\vert^2 \geq
\min_{1\leq j\leq\Ce_m^u} \frac{\lambda_j}{\edfr} \geq\frac
{m^{b_m}}{4\edfr^2 m}\geq
\frac{2}{m},
\]
where the last step holds for $m\geq\exp(128\log(8 \edfr^2)^2)$, as
shown by some algebra.
\end{pf}

For the proof of Lemma~\ref{lem:talalem} below, we need the following
lemma, which can be found in Talagrand \cite{Talagrand96}.
\begin{lem}[(Talagrand's inequality)]
\label{sec:talagrand-1}
Let $T_1, \ldots, T_n$ be independent random variables, and let $\nu^*_n(r)
= (1/n)\sum_{i=1}^n [r(T_i) - \E[r(T_i)] ]$, for $r$
belonging to a countable class $\cR$ of measurable functions. Then,
\begin{eqnarray*}
\E\Bigl[\sup_{r\in\cR} \bigl\vert\nu^*_n(r)\bigr\vert^2 -
6H_2^2 \Bigr]_+ \leq C \biggl(\frac{v}{n}\exp
\bigl(- \bigl(nH_2^2/6v \bigr) \bigr) +
\frac{H_1^2}{n^2} \exp\bigl(-K_2(nH_2/H_1)
\bigr) \biggr)
\end{eqnarray*}
with numerical constants $K_2= (\sqrt{2}-1)/(21\sqrt{2})$ and $C>0$
and with
\[
\sup_{r\in\cR}\Vert r\Vert_\infty\leq H_1, \qquad\E
\Bigl[ \sup_{r\in\cR}\bigl\vert\nu^*_n(r)\bigr\vert\Bigr]\leq
H_2,\qquad\sup_{r\in\cR} \frac{1}{n}\sum
_{i=1}^n \var\bigl(r(T_i) \bigr)\leq
v.
\]
\end{lem}
\begin{lem}\label{lem:talalem}
Let $\delta^*$ and $\Delta^*$ be sequences such
that for all $k\geq1$,
\[
\delta_k^* \geq\sum_{-k\leq j\leq k}
\frac{\omega_j}{\vert[\edf]_j\vert^2} \quad\mbox{and}\quad
\Delta_k^* \geq
\max_{0\leq\vert j\vert\leq k} \frac{\omega_j}{\vert[\edf
]_j\vert^2}
\]
and let $K_2 := (\sqrt{2}-1)/(21\sqrt{2})$. Then, for all $n,k\geq1$,
\begin{eqnarray*}
&&\E\biggl[ \biggl( \supt{k} \bigl\vert\langle t,\tPhi_\nu
\rangle_\hw\bigr\vert^2 - \frac{6 \delta_k^*}{n} \biggr)_+ \biggr]
\\
&&\quad\leq C \biggl\{\frac{\Vert\edf\Vert^2 \Vert\xdf\Vert
^2}{n} \Delta_k^*\exp\biggl(-
\frac{1}{6 \Vert\edf\Vert^2 \Vert\xdf\Vert^2} \bigl(\delta_k^*
/ \Delta_k^* \bigr)
\biggr) + \frac{1}{n^2} \exp( -K_2 \sqrt{n} )
\delta_k^* \biggr\}.
\end{eqnarray*}
\end{lem}

\begin{pf} For $t\in\cS_k$, define the function $r_t := \sum_{k\leq
j\leq k} \hw_j[t]_j\cc{[\edf]}_j^{-1}
\ef_j$. Then it is readily seen that $\langle t,\tPhi_\nu\rangle
_\hw=
\frac{1}{n}\sum_{k=1}^n r_t(Y_k) -
\E[r_t(Y_k)]$. We next compute constants $H_1$, $H_2$, and $v$
verifying the three inequalities required in Lemma~\ref
{sec:talagrand-1}, which then implies the result.

First, consider $H_1$:
\[
\supt{k}\Vert r_t\Vert^2_\infty=
\sup_{y\in\R} \sum_{-k \leq j\leq k} \hw_j
\bigl\vert{ [\edf]}_j\bigr\vert^{-2} \bigl\vert\ef_j(y)
\bigr\vert^2 = \sum_{-k \leq j\leq k} \hw_j
\bigl\vert[\edf]_j\bigr\vert^{-2} \leq\delta_k^*
=:H_1^2.
\]
Next, find $H_2$. Note that
\[
\E\Bigl[\supt{k} \bigl\vert\langle t,\tPhi_\nu\rangle_\hw
\bigr\vert^2 \Bigr] = \frac{1}{n}\sum_{-k\leq j\leq k}
\hw_j\bigl\vert[\edf]_j\bigr\vert^{-2} \var\bigl(
\ef_j(Y_1) \bigr).
\]
Because $\var(\ef_j(Y_1))\leq\E[\vert\ef_j(Y_1) \vert^2]=1$, we
define $
\E[\supt{k} \vert\langle t,\tPhi_\nu\rangle\vert^2] \leq{\delta
_k^*}/{n} =:
H_2^2 $.

Finally, consider $v$. Given $t\in\cB_k$ and a sequence $(z_j)_{j\in
\Z}$, let $\underline{[t]}:=([t]_{-k},\ldots,[t]_k)^T$ and denote by
$D_k(z) := \diag[z_{-k},\ldots,z_k]$ the
corresponding diagonal matrix. Define the Hermitian and positive
semi-definite matrix $A_k := (\cc{[\edf]}_j^{-1} [\edf]_{j'}^{-1}
[\edf]_{j-{j'}} [\xdf]_{j-{j'}}
)_{j,{j'}=-k,\ldots,k}$. Straightforward
algebra shows that
$\supt{k} \var(r_t(Y_1))
\leq\supt{k} \langle A_k D_k(\hw) \vek{[t]},\break D_k(\hw)\vek
{[t]}\rangle_{\C^{2k+1}}$; thus,
\begin{eqnarray*}
\supt{k} \frac{1}{n}\sum_{k=1}^n
\var\bigl(r_t(Y_k) \bigr) &\leq&\supt{k} \bigl\langle
A_k^{1/2}D_k(\hw) \vek{ [t]},A_k^{1/2}D_k(
\hw)\vek{ [t]} \bigr\rangle_{\C^{2k+1}}
\\
&= &\supt{k} \bigl\Vert A_k^{1/2}D_k(\hw)\vek{
[t]}\bigr\Vert_{\C
^{2k+1}}^2 = \bigl\Vert D_k(\sqrt
\hw)A_kD_k(\sqrt\hw)\bigr\Vert_{\C^{2k+1}}.
\end{eqnarray*}
Clearly, we have
$A_k = D_k([\edf]^{-1}) B_k D_k(\cc{[\edf]}^{-1})$,
where $ B_k := ([\edf]_{j-k} [f]_{j-k} )_{j,k=-k,\ldots,k}$.
Consequently,
\begin{eqnarray*}
\supt{k} \frac{1}{n}\sum_{k=1}^n
\var\bigl(r_t(Y_k) \bigr) \leq\bigl\Vert D_k
\bigl(\sqrt{\hw} [\edf]^{-1} \bigr)\bigr\Vert_{\C^{2k+1}}^2
\Vert B_k\Vert_{\C^{2k+1}}.
\end{eqnarray*}
We have that $\Vert D_k(\sqrt{\hw} [\edf]^{-1})\Vert_{\C^{2k+1}}^2
= \max_{0\leq\vert j\vert\leq k} \hw_j\vert[\edf]_j\vert^{-2}
\leq\Delta_k^*$.
It remains to show the boundedness of $\Vert B_k\Vert_{\C^{2k+1}}$. Let
$\ell^2$ be the space of square-summable sequences in $\C$,
and define the operator $B\dvtx\ell^2\to\ell^2$ by $(Bz)_k:= \sum
_{j\in
\Z}[\edf]_{j-k} [f]_{j-k} z_j $, $k\in\Z$. Then it is easily
verified that for any $z\in\ell^2$ with $\Vert z\Vert_{\ell^2}=1$, the
Cauchy--Schwarz inequality yields $ \Vert B z\Vert_{\ell^2}^2 \leq
\Vert\phi\Vert^2 \Vert f\Vert^2$, and thus $\Vert B\Vert_{\ell
^2}^2\leq
\Vert\edf\Vert^2 \Vert\xdf\Vert^2$. Given the orthogonal projection
$\Pi_k$ in $\ell^2$ onto $\cS_k$, the operator $\Pi_k B \Pi_k\dvtx
\cS_k\to\cS_k$ has matrix representation $B_k$ via the isomorphism
$\cS_k\cong\C^{2k+1}$, and hence
$\Vert\Pi_k B \Pi_k\Vert_{\ell^2} = \Vert B_k\Vert_{\C^{2k+1}}$.
Given orthogonal projections with a norm bounded by 1, we conclude that
$\Vert B_k\Vert_{\C^{2k+1}}\leq\Vert B\Vert_{\ell^2}$ for all
$k\in\N
$, which implies that $\supt{k} \frac{1}{n}\sum_{k=1}^n\var
(r_t(Y_k))\leq\Vert\edf\Vert^2 \Vert\xdf\Vert^2 \Delta_k^*=:v$, which
completes the proof.
\end{pf}
\begin{lem}\label{lem:Q}
For every $m\geq1$ and $k\geq0$, we have
\[
\sup_{f\in\Fgr}\E\Bigl[\supt{k}\bigl\vert\langle t,\hPhi_g - \tPhi
_g\rangle_\hw\bigl\vert^2 \Bigr] \leq
C \xdfr\max_{j\in\N} \biggl\{ \frac{\hw_j}{\bw_j}\min\biggl(1,
\frac{1}{m[\edf]_j^2} \biggr)  \biggr\}\leq C d \xdfr\kappa_m(
\xdfw,
\lambda,\hw).
\]
\end{lem}
\begin{pf} First, given that $\xdf\in\cF_\xdfw^\xdfr$, it can be
easily seen that
\begin{eqnarray*}
\E\Bigl[\supt{k} \bigl\vert\langle t,\hPhi_g - \tPhi_g\rangle_\hw
\bigr\vert^2 \Bigr] \leq\xdfr
\sup_{-k\leq j\leq k} \frac{\hw_j}{\xdfw_j} \E\bigl[\vert R_j
\vert^2 \bigr],
\end{eqnarray*}
where $R_j$ is as defined by
\[
R_j := \biggl(\frac{[\edf]_j}{\widehat{[\edf]}_j}\mathbh{1} {
\bigl\{
\bigl\vert\widehat{[\edf]}_j\bigr\vert^2 \geq1/m \bigr\}} -1 \biggr).
\]
In view of the definition \eqref{def:kappam} of
$\kappa_m$, the result follows from $\E[\vert R_j\vert^2] \leq C
\min\{1,
\frac{1}{m\vert[\edf]_j\vert^2} \}$, which can be realized as follows.
Consider the identity
\begin{eqnarray*}
\E\vert R_j\vert^2 = \E\biggl[ \biggl\vert
\frac{[\edf]_j}{\widehat{[\edf]}_j}-1 \biggr\vert^2\mathbh{1} { \bigl\{
\bigl\vert\widehat{[\edf ]}_j\bigr\vert^2 \geq1/m \bigr\}} \biggr] +\P
\bigl[\bigl\vert
\widehat{[\edf]}_j\bigr\vert^2<1/m \bigr] =: R^{\mathrm{I}}_j
+ R^{\mathrm{II}}_j. %
\end{eqnarray*}
Trivially, $R_j^{\mathrm{II}}\leq1$.
If $1\leq4/(m \vert[\edf]_j\vert^2)$, then obviously $R^{\mathrm
{II}}_j\leq
4 \min\{1, \frac{1}{m\vert[\edf]_j\vert^2} \}$.
Otherwise, we have $1/m<\vert[\edf]_j\vert^2/4$ and thus, using
Tchebychev's 
inequality,
\begin{eqnarray*}
R_j^{\mathrm{II}} \leq\P\bigl[\bigl\vert\widehat{[\edf]}_j-[\edf]_j\bigr\vert>
\bigl\vert[\edf]_j\bigr\vert/2 \bigr] \leq
\frac{4 \var(\widehat{[\edf]}_j)}{\vert[\edf]_j\vert^2} \leq4
\min
\biggl\{1, \frac{1}{m\vert[\edf]_j\vert^2} \biggr\},
\end{eqnarray*}
where $\var(\widehat{[\edf]}_j)\leq1/m$ for all $j$.
Now consider $R^\mathrm{I}_j$. We find that
%
%
\begin{equation}
\label{eq:13} R^\mathrm{I}_j = \E\biggl[
\frac{\vert\widehat{[\edf]}_j-[\edf]_j\vert^2}{\vert\widehat
{[\edf ]}_j\vert^2} \mathbh{1} { \bigl\{\bigl\vert\widehat{[\edf]}_j
\bigr\vert^2\geq1/m \bigr\}} \biggr] \leq m\var\bigl( \widehat{[\edf
]}_j\bigr) \leq1.
\end{equation}
On the other hand, using that $\E[\vert\widehat{[\edf]}_j-[\edf
]_j\vert^4]\leq C/m^2$ (cf. Petrov \cite{Petrov1995},
Theorem~2.10), we obtain
\begin{eqnarray*}
R^\mathrm{I}_j &\leq&\E\biggl[ \frac{\vert\widehat{[\edf
]}_j-[\edf
]_j\vert^2}{\vert\widehat{[\edf]}_j\vert^2}
\mathbh{1} { \bigl\{\bigl\vert\widehat{[\edf]}_j\bigr\vert^2\geq1/m
\bigr\}} 2 \biggl\{ \frac{\vert\widehat{[\edf]}_j - [\edf
]_j\vert
^2}{\vert[\edf]_j\vert^2} + \frac{\vert\widehat{[\edf]}_j\vert
^2}{\vert[\edf]_j\vert^2} \biggr\} \biggr]
\\
&\leq&\frac{2 m \E[\vert\widehat{[\edf]}_j-[\edf]_j\vert
^4]}{\vert
[\edf]_j\vert^2} + \frac{2 \var(\widehat{[\edf]}_j)}{\vert[\edf
]_j\vert^2} \leq\frac{2 C}{m \vert[\edf]_j\vert^2} +
\frac{2}{m \vert[\edf]_j\vert^2}.
\end{eqnarray*}
Combining this result with \eqref{eq:13} gives $R^\mathrm{I}_j\leq
2(C+1)\min\{
1, \frac{1}{m\vert[\edf]_j\vert^2} \}$, which
completes the proof.
\end{pf}

\begin{lem}\label{lem:POmegaqc} Under Assumption~\textup{\ref{ass:minreg}},
$\P[\Omega_q^c]\leq C(\edfr) m^{-6} $ for all $m\geq1$.
\end{lem}
\begin{pf} The estimate is obvious for $m< \exp(512\log(3\edfr)^2)=:m_0$.
Consider the complement of $\Omega_q$ given by
\[
\Omega_q^c = \biggl\{\exists0<\vert j\vert\leq
\Ce_m^u \dvt\biggl\vert\frac{[\edf]_j}{\widehat{[\edf]}_j}-1 \biggr\vert>
\frac{1}{2} \vee\bigl\vert\widehat{[\edf]}_j\bigr\vert^2 <1/m
\biggr\}.
\]
Because of Lemma~\ref{lem:NuMu}(iii), for all $m\geq m_0$
and $ 0<\vert j\vert\leq\Ce_m^u$, we have
$\vert[\edf]_j\vert^2\geq2/m$. This yields
\[
\Omega_q^c \subseteq\biggl\{\exists{ 0<}\vert j\vert
\leq\Ce_m^u \dvt\biggl\vert\frac{\widehat{[\edf]}_j}{[\edf]_j} -1
\biggr\vert>
\frac{1}{3} \biggr\}.
\]
By Hoeffding's inequality, for all $ 0<\vert j\vert\leq\Ce_m^u$,
%
%
\begin{equation}
\label{eq:7}\P\bigl[\bigl\vert\widehat{[\edf]}_j/[\edf]_j -1
\bigr\vert>1/3\bigr] \leq2 \exp\biggl(-\frac{m \vert[\edf]_j\vert^2}{72}
\biggr) \leq2\exp\biggl(-
\frac{m \lambda_{\Ce^u_m}}{72\edfr} \biggr),
\end{equation}
which implies the result by virtue of Lemma~\ref{lem:NuMu}(ii).
\end{pf}

\subsubsection*{Fully adaptive estimation}

\begin{pf*}{Proof of Theorem~\ref{thm:upper:unknown}}
We begin the proof by defining the event
$\Omega_{qp}:=\Omega_q\cap\Omega_{p}$, where $\Omega_q$ is given
in \eqref{eq:20} and
%
%
\begin{equation}
\label{app:proofs:gen:2:def:2} \Omega_{p}:= \bigl\{ \bigl(
\Cy^\lambda_n\wedge\Ce^\lambda_m \bigr) \leq
(\hCy_n\wedge\hCe_m) \leq\bigl(\Cy_n^u
\wedge\Ce_m^u \bigr) \bigr\}.
\end{equation}
Observe that on $\Omega_q$, we have $ (1/2)\Delta_k\leq\widehat
\Delta_k\leq(3/2)\Delta_k$ for all $0\leq k\leq
\Ce_m^u$,
and thus $(1/2)[\Delta_k\vee(k+2)]\leq[ \widehat\Delta_k\vee
(k+2)]\leq(3/2)[\Delta_k\vee(k+2)]$, which implies that
\begin{eqnarray*}
&&(1/2) k \Delta_k \biggl(\frac{ \log[ \Delta_k\vee
(k+2)]}{\log(k+2)} \biggr)
\biggl(1-\frac{\log2}{\log
(k+2)}\frac{\log(k+2)}{\log(\Delta_k \vee[k+2])} \biggr)
\\
&&\quad\leq\widehat\delta_k \leq(3/2) k \Delta_k
\biggl( \frac{\log(\Delta_k
\vee[k+2])}{\log( k+2)} \biggr) \biggl(1+\frac{\log3/2}{\log(k+2)
}
\frac{\log(k+2)}{\log(\Delta_k \vee[k+2])} \biggr). %
\end{eqnarray*}
Using $ {\log(\Delta_k
\vee(k+2))}/{\log( k+2)}\geq1$, we conclude from the previous
estimate that
\begin{eqnarray*}
\delta_k/10&\leq& (\log3/2)/(2 \log3)
\delta_k\leq(1/2) \delta_k \bigl[1-(\log2)/\log(k+2)
\bigr] \leq{\widehat\delta_k}
\\
&\leq&(3/2) \delta_k \bigl[1+ (\log3/2)/\log(k+2) \bigr]\leq{ 3
\delta_k}. %
\end{eqnarray*}
Letting $\pen(k):= 60 \delta_k n^{-1}$ and $\hpen(k):= 600
\widehat\delta_k n^{-1}$, it follows that on $\Omega_q$,
\[
\pen(k)\leq\hpen(k)\leq30 \pen(k) \qquad\forall0\leq k\leq\Ce_m^u.
\]
On $\Omega_{qp}=\Omega_q\cap\Omega_{p}$, we have $\whk\leq\Ce
_m^u$. Thus,
%
%
\begin{eqnarray}\label{eq:3}
\bigl( \pen(k\vee\whk) + \hpen(k)-\hpen(\whk) \bigr)\mathbh{1}
{\{
\Omega_{qp}\}} &\leq& \bigl( \pen(k)+\pen(\whk) + \hpen(k)-\hpen
(\whk)
\bigr)\mathbh{1} {\{\Omega_{qp}\}}\hspace*{33pt}
\nonumber
\\[-8pt]
\\[-8pt]
 &\leq&31\pen(k)\qquad\forall0\leq k \leq\Ce_m^u.
\nonumber
\end{eqnarray}
Now consider the decomposition
%
%
\begin{equation}
\label{eq:decomp:qp} \Ex\Vert\hxdf_{\whk}-\xdf\Vert_\omega^2=
\Ex\Vert\hxdf_{\whk}-\xdf\Vert_\omega^2\mathbh{1} {
\{\Omega_{qp}\}} + \Ex\Vert\hxdf_{\whk}-\xdf
\Vert_\omega^2 \mathbh{1} { \bigl\{\Omega_{qp}^c
\bigr\}}.
\end{equation}
We now bound the two terms separately:
\begin{eqnarray*}
\Ex\Vert\hxdf_{\whk}-\xdf\Vert_\omega^2 \mathbh{1}
{\{ \Omega_{qp}\}}&\leq& C \biggl\{ \Vert\xdf- \xdf_k
\Vert_\omega^2 + \edfr\zeta_\edfr\frac{\delta_k^\lambda}{n}
+\xdfr\edfr\kappa_m
\\
&&\hphantom{C \biggl\{}{} + d \zeta_d \frac{\delta_1^\lambda+{\Sigma(\zeta_d^{-1}
\Vert
\edf\Vert^2 \Vert\xdf\Vert^2)}}{n} \biggr\},
\\
\Ex\Vert\hxdf_{\whk}-\xdf\Vert_\omega^2\mathbh{1}
{ \bigl\{\Omega_{p}^c \bigr\}} &\leq& C \biggl(
\frac{\edfr}{\lambda_1} \biggr)^7 \frac{(1+\Vert\xdf\Vert
_\omega^2)}{m}.
\end{eqnarray*}
Consider the first term. Following the proof of \eqref{eq:12}
line by line, it is easily
seen that for $0\leq k\leq(\Cy^\lambda_n\wedge\Ce^\lambda_m)$, we have
\begin{eqnarray*}
&&(1/4) \Vert\hxdf_{\whk}-\xdf\Vert_\omega^2
\mathbh{1} {\{\Omega_{qp}\}}
\\
&&\quad\leq(7/4) (r\hw_k/ \xdfw_k) + 10 \sum
_{j=0}^{\Cy_n^u} \biggl( \supt{j} \bigl\vert\langle t,
\tPhi_\nu\rangle_\omega\bigr\vert^2 - 6
\frac{\delta_{j}}{n} \biggr)_+
\\
&&\qquad{}+ 8\supt{\Cy_n^u\wedge\Ce_m^u}
\bigl\vert\langle t,\hPhi_g - \tPhi_g\rangle_\hw\bigr\vert^2 + \bigl( \pen
(k\vee\whk) + \hpen(k)-\hpen
(\whk)
\bigr) \mathbh{1} {\{\Omega_{qp}\}}
\\
&&\quad\leq(7/4) (r\hw_k/\xdfw_k) + 10\sum
_{j=0}^{\Cy_n^u} \biggl( \supt{j} \bigl\vert\langle t,
\tPhi_\nu\rangle_\omega\bigr\vert^2 - 6
\frac{\delta_{j}}{n} \biggr)_+
\\
&&\qquad{}+ 8\supt{\Cy_n^u\wedge\Ce_m^u}
\bigl\vert\langle t,\hPhi_g - \tPhi_g\rangle_\hw\bigr\vert^2 + 31 \pen(k),
\end{eqnarray*}
where the last inequality follows from \eqref{eq:3}.
The second and third terms are controlled by
Lemmas~\ref{lem:talalem} and~\ref{lem:Q}, respectively (cf. the proof
of \eqref{deco:known:a}). It follows that
%
%
\begin{eqnarray}
\label{deco:unknown:a} \sup_{\xdf\in\Fgr}\sup_{\edf\in\Eld} \E
\Vert
\hxdf_{\hk} -\xdf\Vert^2_\hw\mathbh{1} {\{
\Omega_{qp}\}} &\leq& C \Bigl\{ (\xdfr+ \edfr\zeta_\edfr)
\min_{0\leq k\leq(\Cy
_n^\lambda\wedge
\Ce_m^\lambda)} \bigl[ \max\bigl({\hw_k}/{\xdfw_{k}},
{\delta_k^\lambda}/{n} \bigr) \bigr]
\nonumber
\\[-8pt]
\\[-8pt]
&&\hphantom{C \bigl\{} {}+ \xdfr\edfr\kappa_m + \edfr\Sigma
(\xdfr
\edfr\Edfw\zeta_\edfr) n^{-1} \Bigr\}.
\nonumber
\end{eqnarray}

Consider the second term of \eqref{eq:decomp:qp}.
Following the proof of \eqref{eq:5} and replacing $\Omega_{q}^c$
by $\Omega_{qp}^c$, we obtain
\[
\Ex\Vert\hxdf_{\hk}-\xdf\Vert_\omega^2\mathbh{1}
{ \bigl\{\Omega_{qp}^c \bigr\}} \leq C m^2
\bigl(1+ \Vert f\Vert_\hw^2 \bigr)\P\bigl[
\Omega_{qp}^c \bigr]^{1/2}.
\]

It follows by Lemma~\ref{tech:res2} that for all $m\geq1$,
\[
\sup_{\xdf\in\Fgr}\sup_{\edf\in\Eld} \Ex\Vert\hxdf_{\whk}-
\xdf
\Vert_\omega^2\mathbh{1} { \bigl\{\Omega_{qp}^c
\bigr\}} \leq C(\lambda,\edfr) (1+\xdfr) m^{-1}.
\]
The result of the theorem follows by combining the last estimate
with \eqref{eq:4} and \eqref{deco:unknown:a}.
\end{pf*}

\begin{lem}\label{tech:res2}
Under Assumptions~\textup{\ref{ass:minreg}} and~\textup{\ref
{ass:MplusEins}}, the event
$\Omega_p$ defined in \textup{\eqref{app:proofs:gen:2:def:2}} satisfies
\[
\P\bigl(\Omega_{p}^c \bigr) \leq C(\lambda,\edfr)
m^{-6}\qquad\forall n,m\geq1.
\]
\end{lem}
\begin{pf}Let $\Omega_\mathrm{I}:=\{ (\Cy^\lambda_n\wedge\Ce
^\lambda_m)>(\hCy_n\wedge\hCe_m)
\}$ and $\Omega_{\mathrm{II}}:=\{ (\hCy_n\wedge\hCe_m) > (\Cy
_n^u\wedge
\Ce_m^u)\}$. We then have $\Omega_{p}^c=\Omega_\mathrm{I}\cup
\Omega_{\mathrm{II}}$.
First, consider $\Omega_\mathrm{I} = \{ \hCy_n< (\Cy^\lambda
_n\wedge\Ce^\lambda_m)\} \cup\{
\hCe_m< (\Cy^\lambda_n\wedge\Ce^\lambda_m)\}$. By the definition of
$\Cy^\lambda_n$,
we have that $\min_{1\leq\vert j\vert\leq\Cy^\lambda_n}
\frac{\vert[\edf]_j\vert^{2}}{\vert j\vert\hw_j^+}\geq\frac{4
(\log n)}{n}$, which
implies that
\begin{eqnarray*}
\bigl\{ \hCy_n< \bigl(\Cy^\lambda_n\wedge
\Ce^\lambda_m \bigr) \bigr\}&\subset& \biggl\{ \exists1\leq
\vert j\vert\leq\bigl(\Cy^\lambda_n\wedge
\Ce^\lambda_m \bigr)\dvt\frac{\vert\hfedf_j\vert^2}{\vert
j\vert\hw_j^+}<
\frac{\log n}{n} \biggr\}
\\
&\subset&\bigcup_{1\leq\vert j\vert\leq\Cy^\lambda_n\wedge\Ce
^\lambda_m} \biggl\{ \frac{\vert\hfedf_j\vert}{\vert[\edf]
_j\vert}
\leq1/2 \biggr\} \subset\bigcup_{1\leq\vert j\vert\leq\Cy
^\lambda
_n\wedge
\Ce^\lambda_m} \biggl\{ \biggl\vert
\frac{\hfedf_j}{[\edf]_j}-1 \biggr\vert\geq1/2 \biggr\}.
\end{eqnarray*}
From $\min_{1\leq\vert j\vert\leq\Ce^\lambda_m}
\vert[\edf]_j\vert^{2}\geq4 m^{-1+ b_m}$, it follows in the same
way that
\begin{eqnarray*}
\bigl\{ \hCe_m< \bigl(\Cy^\lambda_n\wedge
\Ce^\lambda_m \bigr) \bigr\}\subset\bigcup
_{1\leq\vert j\vert\leq\Cy^\lambda_n\wedge
\Ce^\lambda_m} \biggl\{ \biggl\vert\frac{\hfedf_j}{[\edf]_j}-1 \biggr\vert
\geq1/2 \biggr
\}.
\end{eqnarray*}
Therefore, $\Omega_\mathrm{I} \subset\bigcup_{1\leq\vert j\vert
\leq
\Ce_m^u} \{ \vert\hfedf_j/[\edf]_j-1\vert\geq1/2 \}$, because
$\Ce
_m^\lambda
\leq\Ce_m^u$. Thus, applying \mbox{Hoeffding's} inequality and
Lemma~\ref{lem:NuMu}(ii) as in \eqref{eq:7} yields
%
%
\begin{equation}
\label{tech:e1} \P[\Omega_\mathrm{I}] \leq\sum
_{1\leq\vert j\vert\leq\Ce_m^u} 2 \exp\biggl(-\frac{m \vert
[\edf]_j\vert^2}{72} \biggr) \leq C(
\edfr) m^{-6}.
\end{equation}
Consider $\Omega_{\mathrm{II}} = \{ \hCy_n> (\Cy_n^u\wedge\Ce
_m^u)\} \cap
\{ \hCe_m> (\Cy_n^u\wedge\Ce_m^u)\}$.
In the case where $(\Cy_n^u\wedge\Ce_m^u)= \Cy_n^u$, use $\frac
{\log n}{4n} \geq
\max_{\vert j\vert\geq\Cy_n^u +1} \frac{\vert[\edf]_j\vert
^2}{\vert j\vert\hw_j^+}$, such that
\begin{eqnarray*}
\Omega_{\mathrm{II}}&\subset& \bigl\{\hCy_n> \Cy_n^u
\bigr\}\subset\biggl\{ \forall1\leq\vert j\vert\leq\Cy_n^u
+ 1\dvt\frac{\vert\hfedf_j\vert^2}{\vert j\vert\hw_j^+}\geq
\frac{\log n}{n} \biggr\}
\\
&\subset& \biggl\{ \frac{\vert\hfedf_{\Cy_n^u+1}\vert}{\vert
[\edf]_{\Cy
_n^u+1}\vert}\geq2 \biggr\}\subset\bigl\{ \bigl\vert
\hfedf_{\Cy_n^u+1}/[\edf]_{\Cy_n^u+1}-1\bigr\vert\geq1 \bigr\}.
\end{eqnarray*}
In the case where $(\Cy_n^u\wedge\Ce_m^u)= \Ce_m^u$, it follows
analogously from $m^{-1+b_m} \geq
4\times\max_{\vert j\vert\geq\Ce_m^u+1} \vert[\edf]_j\vert^2$ that
\begin{eqnarray*}
\Omega_{\mathrm{II}}\subset\bigl\{ \hCe_m> \Ce_m^u
\bigr\} \subset\bigl\{ \bigl\vert\hfedf_{\Ce_m^u+1}/[\edf]_{\Ce
_m^u+1}-1\bigr\vert\geq1
\bigr\}.
\end{eqnarray*}
Therefore, we have $\Omega_{\mathrm{II}} \subset\{
\vert\hfedf_{(\Cy_n^u\wedge\Ce_m^u)+1}/[\edf]_{(\Cy_n^u\wedge
\Ce
_m^u)+1}-1\vert\geq
1 \}$. Applying Hoeffding's
inequality as in \eqref{eq:7} and using
Assumption~\ref{ass:MplusEins}, we obtain, for all $m\geq1$,
%
%
\begin{equation}
\label{tech:e2} \P[\Omega_{\mathrm{II}}] \leq2 \exp\biggl(-
\frac{m \vert[\edf]_{\Ce_m^u+1}\vert^2}{72} \biggr) \leq
C(\lambda
,\edfr) m^{-7}.
\end{equation}
Combining \eqref{tech:e1} and \eqref{tech:e2} implies the result.
\end{pf}
%

\subsection*{Illustration: Estimation of derivatives}
\label{sec:illustration}

\begin{pf*}{Proof of Proposition~\ref{prop:ex-cont-ada-ukn}} In light
of the
proof of Proposition~\ref{coro:ex:lower}, we apply
Theorem~\ref{thm:upper:unknown}, where in both cases we need only
check the additional Assumption~\ref{ass:MplusEins}. The result then follows
by an evaluation of the upper bound.

{[o-o]} It is easily seen that $(m \lambda_{\Ce
_{m}^u+1})^{-1}\log
m= o(1)$ as $m\to\infty$. Thus, Assumption~\ref{ass:MplusEins}
is satisfied in this case. Because $\kstar\sim n^{1/(2a+2p+1)}$, we have
$\kstar\lesssim\Cy_n^\lambda$. Thus, the upper bound is
%
%
\begin{equation}
\label{pr:pro:ada:kno:poly:e1} \bigl(\kstar\wedge\Mmn^\lambda
\bigr)^{-2(p-s)} + m_n^{-(1\wedge((p-s)/a))}.
\end{equation}
We consider two cases. First, let $p-s>a$. Suppose that
$n^{2(p-s)/(2p+2a+1)} = \mathrm{O}(m_n)$; then,
\[
\frac{\kstar}{\Mmn^\lambda} \sim\frac{n^{1/(2a+2p+1)}}{
(m_n^{1-b_{m_n}} )^{(1/2a)}} = \frac{n^{1/(2a+2p+1)}}{m_n^{1/2(p-s)}}
\bigl(m_n^{-a+(p-s)(1-b_{m_n})} \bigr)^{{1}/{(2(p-s)a)}} = \mathrm{o}(1).
\]
This means that $\kstar\lesssim\Mmn^\lambda$, so the resulting
upper bound
is $(\kstar)^{-2(p-s)} + m_n^{-1} \lesssim(\kstar)^{-2(p-s)}$.
Suppose now that $m_n = \mathrm{o}(n^{2(p-s)/(2p+2a+1)})$. If in addition
$\kstar= \mathrm{O}( \Mmn^\lambda)$, then the first summand in
\eqref{pr:pro:ada:kno:poly:e1} reduces to $(\kstar)^{-2(p-s)}$ and
thus the upper bound is $m_n^{-1}$. On the other hand, if
$\Mmn^\lambda/\kstar= \mathrm{o}(1)$, then the first term is
$(\Mmn^\lambda)^{-2(p-s)}\sim(m_n^{-a+(p-s)(1-b_{m_n})})^{1/a}
m_n^{-1}\lesssim m_n^{-1}$,
because $p-s>a$. Combining both cases, we obtain the result in the case where
$p-s>a$.

Now assume that $p-s\leq a$. First, suppose that $\kstar=
\mathrm{O}(\Mmn^\lambda)$. Then the first summand in
\eqref{pr:pro:ada:kno:poly:e1} reduces to $(\kstar)^{-2(p-s)}$, and,
moreover, it follows that $n^{2a/(2p+2a+1)} = \mathrm{O}(m_n)$.
Therefore, the
upper bound is
$(\kstar)^{-2(p-s)}$. Now consider $\Mmn^\lambda= \mathrm{o}(\kstar
)$. Then
\eqref{pr:pro:ada:kno:poly:e1} can be rewritten as
$(m_n^{1-b_{m_n}})^{-(p-s)/a} + m_n^{-(p-s)/a}$, which results in the rate
$(m_n^{1-b_m})^{-(p-s)/a}$. Combining both cases gives the
result. More precisely, $m_n = \mathrm{o}(n^{2a/(2p+2a+1)})$ implies
$\Mmn^\lambda=
\mathrm{o}(\kstar)$. In contrast, in the case where $n^{2a/(2p+2a+1)}
= \mathrm{O}(m_n)$,
if $\kstar/\Mmn^\lambda= \mathrm{O}(1)$, then the rate is $(\kstar)^{-2p}$,
whereas if
$\Mmn^\lambda/\kstar= \mathrm{o}(1)$, then the rate is
$(m_n^{1-b_m})^{-(p-s)/a}$.

{[s-o]} As in case [o-o], Assumption~\ref{ass:MplusEins}
is satisfied. Recall that $\kstar\sim(\log n)^{1/(2p)}$. If~$n(\log
n)^{-(2a+2s+1)/(2p)}= \mathrm{O}(m_n)$, then $\kstar\lesssim\Ce
_{m_n}^\lambda$
and $m_n^{-1}\lesssim\psi_{n,m_n}^\diamond\sim\break n^{-1}(\log
n)^{(2a+2s+1)/(2p)}$.
In the
opposite case, we have $\psi_{n,m_n}^\diamond\lesssim m_n^{-1}$,
which proves the result.

{[o-s]} To verify that Assumption~\ref{ass:MplusEins} is
satisfied in this setting, we can proceed as follows. Define the
sequence $\tCe^u$ exactly as $\Ce^u$, but replacing $b_m$ by $a_m =
b_m^{2^k}$. Then $\tCe^u$ satisfies assertion
Lemma~\ref{lem:NuMu}(ii), the proof being similar to that for
$\Ce^u$. In contrast, we can show that $\tCe^u_m -
\Ce^u_m\to\infty$ as $m\to\infty$, which amounts to showing
Assumption~\ref{ass:MplusEins}.

We have $\kstar\asymp(\log n)^{1/2a}$. The upper bound becomes
$(\kstar\wedge\Mmn^\lambda)^{-2(p-s)} +\break (\log m_n)^{-(p-s)/a}
\sim(\kstar\wedge\Mmn^\lambda)^{-2(p-s)}$. Distinguishing
$\kstar\lesssim\Ce^\lambda_m$ and the opposite case shows the result.
\end{pf*}
\end{appendix}

\section*{Acknowledgements}

This work was supported by the IAP research network no. P6/03 of the
Belgian Government (Belgian Science Policy) and by the ``Fonds
Sp\'eciaux de Recherche'' from the Universit\'e catholique de Louvain.

We are grateful to two referees for raising a number of points which
helped us clarify several important issues. Finally, we thank the
Associate Editor for constructive criticism and clear guidelines.


%

%

\printhistory

\end{document}